\documentclass[a4paper,11pt]{article}

\pdfoutput=1

\usepackage{graphicx}

\usepackage{color}
\usepackage{authblk}
\usepackage{amsmath,amscd,amsthm,amssymb}
\numberwithin{equation}{section}

\usepackage{algorithm2e}
\usepackage{caption}
\usepackage{sub caption, capt-of}
\usepackage{mathrsfs}

\graphicspath{ {./Figures/} }

\definecolor{grey}{gray}{.35}

\newtheorem{theorem}{Theorem}[section]
\newtheorem{Definition}[theorem]{Definition}

\newtheorem{Remark}[theorem]{Remark}
\newtheorem{Algorithm}[theorem]{Algorithm}

\newcounter{listitem}

\newtheoremstyle{thmstl}{\topsep}{\topsep}
 {}%         Body font
 {}%         Indent amount (empty = no indent, \parindent = para indent)
 {\bfseries}% Thm head font
 {}%        Punctuation after thm head
 {2ex}%     Space after thm head (\newline = linebreak)
 {\thmname{#1}\thmnumber{ #2}}%         Thm head spec
\theoremstyle{thmstl}

\newcounter{hyp}
\newtheorem{hypothesis}{}

\newcounter{assump}
\newtheorem{assumption}{}

\def\R{{\mathbb R}}
\def\Z{{\mathbb Z}}

\makeatletter \@addtoreset{equation}{section}

\makeatother

\def\qed{\hbox{\hskip 6pt\vrule width7pt height7pt depth1pt
\hskip1pt}\bigskip}

\begin{document}{}

\title{Approximation of reachable sets\\ using optimal control and support vector machines}

\author[1]{Martin Rasmussen\thanks{m.rasmussen@imperial.ac.uk}}
\author[2]{Janosch Rieger\thanks{rieger@math.uni-frankfurt.de}}
\author[1]{Kevin N.~Webster\thanks{kevin.webster@imperial.ac.uk}}
\affil[1]{\small Department of Mathematics, Imperial College London, 180 Queen's Gate, London SW7 2AZ, UK}
\affil[2]{\small FB Mathematik, Goethe-Universit\"{a}t Frankfurt, Postfach 11 19 32, 60054 Frankfurt am Main, Germany}
\maketitle
\vspace{-1cm}
\begin{abstract}
We propose and discuss a new computational method for the numerical approximation of reachable sets for nonlinear control systems.
It is based on the support vector machine algorithm and represents the set approximation as a sublevel set of a function chosen
in a reproducing kernel Hilbert space.
In some sense, the method can be considered as an extension to the optimal control algorithm approach recently developed by Baier,
Gerdts and Xausa.
The convergence of the method is illustrated numerically for several examples.
\end{abstract}

%\textbf{Keywords:} Reachable set, Control system, Optimal control, Support vector machine

\section{Introduction}\label{sec:introduction}

The numerical computation of reachable sets is a crucial topic in nonlinear control theory and the quantification of
deterministic uncertainty in dynamical systems.
Collision avoidance of manned and unmanned vehicles is one particular application that currently attracts a lot of attention
(see e.g.\ \cite{Gerdts_12} and the references therein).
Standard techniques such as the set-valued Euler method \cite{Dontchev_89_1,Beyn_07_1} evolve a grid-based approximation
of the reachable set along the relevant time interval.
They are typically very slow, because there is a high degree of redundancy in the computations they carry out.

Recently, a version of the set-valued Euler method was presented in \cite{Rieger_Unpub_2} that tracks the boundaries
of the reachable sets and uses only the boundaries of the right-hand side of the differential inclusion.
With this approach, the complexity of the Euler scheme is reduced drastically in the low-dimensional setting, but only
marginally in higher dimensions.

The DFOG optimal control algorithm \cite{Baier_13_1}, which will be discussed in more detail in Section \ref{section:dfog},
is another recent attempt to reduce the proportion of irrelevant computations.
Every point of a grid in the relevant region of the phase space is projected to the reachable set by solving a Mayer problem.
From this data, one can derive -- at least theoretically -- an accurate description of the reachable set.
In contrast to traditional methods, there is no guarantee that the numerical optimisation routine
finds a global minimum, and therefore, the algorithm is, strictly speaking, unstable.
Numerical studies, however, support the usefulness of this method.

In this paper, we propose a new approach to the calculation and representation of a reachable set approximation,
motivated as an extension to the DFOG algorithm.
The extension consists of using the results of these optimal control problems to search for a function in a particular
function space, so that the reachable set is represented as a sublevel set of this function.
The function space under consideration is a reproducing kernel Hilbert space (RKHS), and the algorithm to search for this function
is a modified support vector machine (SVM) algorithm.

Our algorithm has the advantage that it is robust to a small number of errors made by the optimisation routines from the DFOG method.
In addition, the function used for the reachable set approximation has a sparse representation in terms of the optimal control
results, and the algorithm focuses on information provided by points that are close to the boundary of the reachable set.

\section{Reachable sets and known techniques for their approximation}

In the following, we give a condensed overview over basic properties of reachable sets (see Section \ref{sec:reachable}),
the currently most common numerical methods for approximating them (see Section \ref{sec:Runge:Kutta}) and
the DFOG method (see Section \ref{section:dfog}), which is the basis of our new method.

We recall some standard definitions with regard to set representations.
\begin{Definition}
  Let $A,B\subset\mathbb{R}^d$ be compact sets, and $x\in\mathbb{R}^d$. The distance of a point $x$ to the set $A$ is defined by
  \[\mbox{\normalfont dist}(x,A):= \mbox{\normalfont inf}_{a\in A} \|x-a\|\,.\]
  For any $r>0$, the $r$-neighbourhood of $A$ is the set
  \[B(A,r):=\{z\in\R^d: \mbox{\normalfont dist}(z,A)\le r\}.\]
  The projection of $x$ to $A$ is the set of points in $A$ that realise the infimum distance to $x$, i.e.
  \[\operatorname{Proj}(x,A) :=\{a\in A:\|x-a\| = \mbox{\normalfont dist}(x,A)\}\,.\]
  The Hausdorff semi-distance between sets $A$ and $B$ is given by
  \[\mbox{\normalfont d}(A,B):= \mbox{\normalfont sup}_{a\in A}\mbox{\normalfont dist}(a,B),\]
  and the Hausdorff distance between $A$ and $B$ is given by
  \[\mbox{\normalfont d}_H(A,B):= \mbox{\normalfont max}\{\mbox{\normalfont d}(A,B), \mbox{\normalfont d}(B,A)\}\,.\]
\end{Definition}
Throughout this paper, the symbol $\|\cdot\|$ denotes the Euclidean norm.
The symbols $\|\cdot\|_\infty$, $\mbox{\normalfont dist}_\infty(x,A)$ etc.\ denote the corresponding concepts
based on the maximum norm.

\subsection{Reachable sets} \label{sec:reachable}

Let $U$ be a nonempty convex and compact subset of $\mathbb{R}^d$ and
\begin{displaymath}
  \mathcal{U}:=\big\{u\in L^\infty([t_0,T],\mathbb{R}^d):u(t)\in U \text{ for almost all } t\in [t_0,T]\big\}
\end{displaymath}
for fixed times $t_0 < T$. We consider the nonlinear control problem
\begin{subequations} \label{gesamt}
 \begin{align}
\dot{x}(t) &=  g(t,x(t),u(t))\,,\quad 	u \in  \mathcal{U}, \label{eqn:controlproblem1}\\
x(t_0) &= x_0, \label{eqn:controlproblem2}
\end{align}
\end{subequations}
for some $x_0\in\mathbb{R}^d$, where \eqref{eqn:controlproblem1} holds for almost every $t\in[t_0,T]$ and $x(\cdot)\in W^{1,\infty}([t_0,T],\mathbb{R}^d)$ is absolutely continuous. We are interested in the \textit{reachable set} at time $T$, given by
\begin{displaymath}
\mathcal{R}(T,t_0,x_0) := \{x(T) : x(\cdot) \text{ solves \eqref{gesamt}}\} \,.
\end{displaymath}
Problem \eqref{gesamt} is equivalent to the differential inclusion
\begin{subequations} \label{eqn:controlprobleminclusion}
\begin{align}
\dot{x}(t) & \in  G(t,x(t))\label{eqn:controlprobleminclusion1}\\
x(t_0) & =  x_0, \label{eqn:controlprobleminclusion2}
\end{align}
\end{subequations}
with \eqref{eqn:controlprobleminclusion1} valid for almost all $t\in [t_0,T]$, and $G(t,x) := \bigcup_{u\in U} \{g(t, x, u)\}$.

Reachable sets of nonlinear control systems, or, equivalently, nonlinear differential inclusions, are, in general, nonconvex.
It is, however, well-known, that they enjoy several favourable properties under mild assumptions imposed on the right-hand side
(see e.g.~\cite[Corollary 7.1]{Deimling_92_1}):

\begin{theorem}
  Let $G:[t_0,T]\times\R^d\rightrightarrows\R^d$ have closed and convex images, and assume that
  \begin{itemize}
  \item [a)] the mapping $t\mapsto G(t,x)$ is measurable for all $x\in\R^d$,
  \item [b)] the mapping $x\mapsto G(t,x)$ is upper semicontinuous for all $t\in[t_0,T]$,
  \item [c)] there exists $c\in L^1([t_0,T])$ such that $\|G(t,x)\|\le c(t)(1+\|x\|)$ for all $t\in[t_0,T]$ and $x\in\R^d$.
  \end{itemize}
  Then the mapping $x_0\mapsto\mathcal{R}(T,t_0,x_0)$ is upper semicontinuous, and the reachable set $\mathcal{R}(T,t_0,x_0)$
  is nonempty and compact for all $x_0\in\R^d$.
\end{theorem}

\subsection{Runge-Kutta methods} \label{sec:Runge:Kutta}

Reachable sets may be approximated numerically using set-valued Runge-Kutta methods.
Given a time discretisation $t_n = t_0 + nh$ with $n\in\{0,\dots, N\}$ and $h=T/N$, the iterations
\begin{equation} \label{svrkm}
x_{n+1} \in x_n + h\Phi(t_n,x_n),\quad n=0,\ldots,N,
\end{equation}
with initial value $x_0$ from \eqref{eqn:controlprobleminclusion2} and
\begin{eqnarray*}
\mathcal{R}_h(t_{n+1},t_0,x_0) &=& \bigcup_{x\in\mathcal{R}_h(t_n,t_0,x_0)}\,x+h\Phi(t_n,x),\quad n=0,\ldots,N-1	 \label{eqn:multivaluedRK1}\\
\mathcal{R}_h(t_0,t_0,x_0) &=& \{x_0\},	\label{eqn:multivaluedRK2}
\end{eqnarray*}
with a suitable increment function $\Phi: [t_0,T]\times\R^d\rightrightarrows\R^d$ define the trajectories and the reachable sets
of the numerical scheme.
%, which approximate trajectories and reachable sets at the corresponding times of the original differential
%inclusion.
The simplest example of such a Runge-Kutta scheme is the set-valued Euler method with increment function $\Phi(t,x)=G(t,x)$,
which has been studied in \cite{Dontchev_89_1,Dontchev_00_1} and several other contributions.
The central result for our purposes is published in \cite{Dontchev_89_1}.
\begin{theorem} \label{semi:discrete:Euler}
Let $G:\R^d\rightrightarrows\R^d$ be Lipschitz continuous w.r.t.\ $\mbox{\normalfont d}_H$ with convex and compact values.
Then the reachable sets of the Euler scheme satisfy
\begin{equation*}
d_H(\mathcal{R}(T, t_0, x_0), \mathcal{R}_h(T, t_0, x_0)) \le \text{const}\cdot h.
\end{equation*}
\end{theorem}

In practice, to compute approximations of the reachable set it is necessary to spatially discretise these schemes.
The most natural approach is to introduce a grid $\rho\mathbb{Z}^d$ in the phase space $\R^d$ with grid size $\rho>0$,
and to define a fully discretized scheme with trajectories
\[x_{n+1} \in B_\infty(x_n + h\Phi(t_n,x_n),\rho/2)\cap\rho\Z^d,\quad n=0,\ldots,N,\]
yielding discrete reachable sets
\begin{subequations}
\begin{eqnarray}
\mathcal{R}_{h,\rho}(t_{n+1},t_0,x_0)
&=& \bigcup_{x\in\mathcal{R}_{h,\rho}(t_n,t_0,x_0)}\,B_\infty(x+h\Phi(t_n,x),\rho/2)\cap\rho\Z^d,\quad n=0,\ldots,N-1, \label{union}\\
\mathcal{R}_{h,\rho}(t_0,t_0,x_0) &=& B_\infty(x_0,\rho/2)\cap\rho\Z^d.
\end{eqnarray}
\end{subequations}
The blowup of the images is necessary in order to obtain subsets of the grid that are
close to the original Euler sets in Hausdorff distance.
Error estimates corresponding to spatial discretisation have been studied in \cite{Beyn_07_1}.
Some results are subsumed in the following statement.
\begin{theorem} \label{fully:discrete:thm}
Let $G:\R^d\rightrightarrows\R^d$ be Lipschitz continuous w.r.t.\ $\mbox{\normalfont d}_H$
with convex and compact values.
Then the reachable sets of the fully discrete Euler scheme satisfy
\[\mbox{\normalfont d}_H(\mathcal{R}(T,t_0,x_0),\mathcal{R}_{h,\rho}(T,t_0,x_0))
\le \text{const}\cdot(\rho+h+\rho/h).\]
\end{theorem}
The error term $\text{const}\cdot\rho/h$ forces the user to choose a very fine spatial discretisation,
causing a high computational complexity, which is worsened by the high level of redundancy incurred by
computing parts of the reachable set over and over again in the union \eqref{union}.
Higher-order Runge--Kutta methods are practically infeasible when directly transferred from
ordinary differential equations to inclusions, and hence do not seem to be a cure for the complexity problem,
because the computational costs of a successive evaluation of multifunctions undo all positive effects of the higher
order time-discretisation.

Therefore, Euler's method has been the main focus of study in this area.
A method that is based on the Euler scheme and tracks the boundary instead of the complete reachable set
for reducing the computational cost has been studied in \cite{Rieger_Unpub_2}.
Variations of the implicit Euler scheme that are superior to the explicit Euler scheme when applied to stiff
differential inclusions have been analysed in \cite{Beyn_10_1} and \cite{Rieger_14_1}.

\subsection{The DFOG method}\label{section:dfog}

In this section we review a different approach that has been recently proposed by Baier, Gerdts and Xausa \cite{Baier_13_1},
and which they call the DFOG method (short for \emph{distance fields on grids}).
Exploratory work has been published earlier in \cite{BaiBusChaGer07} and \cite{Baier_09_1}.
This method exploits the representation
\begin{equation*}
  A = \mathbb{R}^d \ \backslash\  \bigcup_{x\in \mathbb{R}^d} \mbox{\normalfont int } B(x,\operatorname{dist}(x,A))\,.
\end{equation*}
of a closed set $A$ as the complement of the union of all open balls contained in $A^c$.
% \cite{Clarke_98_1}:
Note that $x\in A$ implies $\operatorname{int} B(x,\mbox{\normalfont dist}(x,A)) = \emptyset$.

\begin{Algorithm}[DFOG method] \label{DFOGmethod}
Let $x_0\in\mathbb{R}^d$ and $t_0<T$ be as in \eqref{eqn:controlproblem2}.
Let $\rho>0$, let $\Omega\subset\mathbb{R}^d$ be a bounded region with $B_\infty(\mathcal{R}(T,t_0,x_0),\rho/2)\subset\Omega$,
and define a grid $\tilde\Omega:=\Omega\cap\rho\mathbb{Z}^d$.

For every $z\in \tilde\Omega$, solve the optimisation problem
\begin{equation*}
\operatorname{Minimise} \frac{1}{2}\|x_N - z\|^2
\end{equation*}
over  trajectories $(x_n)_{n=0}^N$ of the Runge-Kutta scheme \eqref{svrkm}.
Let $(x_n^*(z))_{n=0}^N$ be a solution to this problem, and denote $x^*(z):=x_N^*(z)$
and $\theta(z) := \frac{1}{2}\|x_N^*(z) - z\|^2$.
Then the set
\begin{equation} \label{eqn:DFOGrepresentation}
  \mathcal{R}_{DFOG}(T,t_0,x_0):= \Omega \ \backslash \  \bigcup_{z\in\tilde\Omega}
  \mbox{\normalfont int } B(z,\sqrt{2\theta(z)})\,.
\end{equation}
is an approximation to the reachable set $\mathcal{R}(T,t_0,x_0)$.
\end{Algorithm}

Algorithm \ref{DFOGmethod} requires the solution of gobal optimisation problems in a very high-dimensional state space,
to which standard tools for local optimisation such as the SQP method may be applied (see for example \cite{Nocedal_99_1}).
Due to the high state-space dimension, it is impossible to use global optimisation routines.
Therefore, one potential pitfall is the existence of local minima leading to incorrect results for $x^*(z)$ and $\theta(z)$.
This means that the representation \eqref{eqn:DFOGrepresentation} may potentially cut away large parts of the reachable set.

One heuristic solution to this problem offered in \cite{Baier_13_1} is `ball-checking':
For any $z,z'\in\tilde\Omega$, check whether the computed $x^*(z)$ and $\theta(z')$ satisfy
$x^*(z)\in B(z',\sqrt{2\theta(z')})$.
In that case, the optimisation routine failed to compute $\theta(z')$ properly, so that
the ball $B(z',\sqrt{2\theta(z')})$ is incorrect and must be ignored.
However, this strategy does not necessarily detect erroneous results.

The following statement is a consequence of Theorem \ref{semi:discrete:Euler}
and the representation \eqref{eqn:DFOGrepresentation}.
\begin{theorem} \label{dfog:estimate}
Let $G:\R^d\rightrightarrows\R^d$ be Lipschitz continuous w.r.t.\ $\mbox{\normalfont d}_H$ with convex and compact values.
If the optimisation problems in Algorithm \ref{DFOGmethod} are solved correctly, then
\begin{equation*}%\label{eqn:DFOGestimate_rhoh}
\mbox{\normalfont d}_H(\mathcal{R}(T,t_0,x_0),\mathcal{R}_{DFOG}(T,t_0,x_0)) \le \text{const}\cdot(h+\rho)\,.
\end{equation*}
\end{theorem}
In contrast to Theorem \ref{fully:discrete:thm}, the error estimate does not contain the critical term $\text{const}\cdot\rho/h$,
which indicates that the DFOG method should be substantially faster than the Euler scheme.
Nevertheless, it is difficult to compare the performance of both schemes on a formal level, because their design and their
behaviour are too different.

The DFOG method can be accelerated using the method of maximal gains published in \cite{Rieger_13_1},
which does not use a spatial a grid, but chooses optimal test points with respect to all information
that is available at runtime.

\section{SVM Algorithm}

In this section, we develop step by step the algorithm the present paper is concerned with.
As the DFOG method, it requires the results of the optimal control routine given in Algorithm \ref{DFOGmethod},
but it uses the well-known SVM algorithm from machine learning to represent the approximation of the reachable set as a sublevel set
of a smooth function chosen from a reproducing kernel Hilbert space (RKHS). % instead of the complement of balls.
This representation provides a smoother boundary for the reachable set approximation and %provides
some robustness against a small number of errors corresponding to the global optimisation routine finding local minima.

In Section \ref{sec:rkhs}, we give a brief introduction to reproducing kernel Hilbert spaces,
and in Section \ref{sec:SVM}, we adapt the classical SVM algorithm to our particular problem.
In Section \ref{sec:incremental}, we discuss the possibility of an adaptive enlargement of the
dataset and computational implications.
A brief comment on the ball-checking procedure in the SVM context is given in Section \ref{sec:ball:check}.

\subsection{Reproducing Kernel Hilbert Space} \label{sec:rkhs}

We will define our RKHS in terms of a Mercer kernel.
\begin{Definition}
Let $X\subset\mathbb{R}^d$.
A Mercer kernel is a function $K: X \times X \rightarrow\mathbb{R}$ satisfying
\begin{itemize}
\item [a)] $K(x,x')=K(x',x)$ for all $x,x'\in X$,
\item [b)] $\sum^n_{i,j=1} c_ic_jK(x_i,x_j)\ge0$ for any $n\in\mathbb{N}$, any $c_1,\ldots,c_n\in\mathbb{R}$
and any $x_1,\ldots,x_n\in X$.
\end{itemize}
\end{Definition}

Typical examples of kernel functions include the Gaussian kernel
\begin{equation}
K_G(x,y) = \textrm{exp}\big(-\textstyle\frac{1}{\sigma}\|x-y\|^2\big)\qquad \text {with } \sigma >0\,,	\label{eqn:Gaussiankernel}
\end{equation}
and the \textit{degree-}$p$ polynomial kernel
\begin{equation*} %\label{eqn:polynomialkernel}
K_P(x,y) = (x^T y +\tau)^p \qquad \text{with } \tau\ge 0\,.
\end{equation*}

Given a kernel function, we define $K_x:=K(x,\cdot)$. The following theorem states how a Mercer kernel uniquely defines
a reproducing kernel Hilbert space \cite{Aronszajan_50_1}:

\begin{theorem}[Moore--Aronszajn]	\label{thm:Moore-Aronszajn}
Given a Mercer kernel $K$, there exists a unique Hilbert space $\mathcal{H}_K$ of functions on $X$ with associated inner product
$\langle\cdot,\cdot \rangle_K$ satisfying the following conditions:
\begin{itemize}
\item[(i)] $K_x\in\mathcal{H}_K$ for all $x\in X$,
\item[(ii)] $\mbox{\normalfont span}\{K_x :   x\in X\}$ is dense in $\mathcal{H}_K$,
\item[(iii)] $f(x) = \langle K_x,f\rangle_K$ for all $f\in\mathcal{H}_K$ and all $x\in X$.
\end{itemize}
\end{theorem}
The inner product in the RKHS is defined by $\langle K_x,K_y\rangle_K = K(x,y)$ and extending linearly. $\mathcal{H}_K$ is then taken
as the completion of the linear span of $\{K_x :  x\in X\}$ with respect to this inner product. The third property in Theorem
\ref{thm:Moore-Aronszajn} is the \textit{reproducing property}.

In this paper, we choose to work with the RKHS corresponding to the Gaussian kernel \eqref{eqn:Gaussiankernel}, although the
algorithm is also viable with other choices of kernels.
The Gaussian RKHS has been well studied and is a very rich function space to work in, which is illustrated
by the following result from \cite{Steinwart_01_1}.

\begin{theorem}
Let $X\subset\mathbb{R}^d$ be compact. Then the Gaussian RKHS $\mathcal{H}_K$ on $X$ is dense in the space $C(X)$
of continuous functions on $X$.
\end{theorem}

For a detailed coverage of the Gaussian RKHS, we refer to \cite{Minh_10_1,Steinwart_06_1}.
In practice, we will be working with finite dimensional RKHS, and in particular, the set $X$ will be chosen according to the grid
points and results from the DFOG optimal control method.

\subsection{Support Vector Machine} \label{sec:SVM}

Support vector machines (SVMs) are well-known supervised learning algorithms frequently used for classification problems,
a common task in machine learning problems.
The soft-margin SVM algorithm was first proposed by Cortes and Vapnik \cite{Cortes_95_1} and is now a popular choice
for machine learning problems.
Applications of the SVM algorithm include handwriting recognition, image classification and text categorisation
\cite{DeCoste_02_1, Joachims_97_1, Schoelkopf_01_1, Tong_01_1}.

We apply the SVM algorithm to a labelled training set $\mathcal{D}$, which contains all relevant information encoded in
the output data $(z,\theta(z),x^*(z))_{z\in\tilde\Omega}$ of the DFOG method, in order to recognise the shape of the
reachable set. The training set $\mathcal{D}$ and index sets $\mathcal{I}$, $\mathcal{E}$ and $\mathcal{B}$
that partition $\mathcal{D}$ into interior, exterior and boundary points of the reachable set, respectively,
are constructed as follows.

\begin{Algorithm}[Labelling]  \label{construct:training:set}
First run Algorithm \ref{DFOGmethod} to obtain the data $(z,\theta(z),x^*(z))_{z\in\tilde\Omega}$. Set $\mathcal{I},\mathcal{E},\mathcal{B}:=\emptyset$ and $m:=0$. Fix $\epsilon>0$.\\
{\normalfont for} $z\in\tilde\Omega$\\
\mbox{}\hspace{3ex} {\normalfont if} $\theta(z)\le\epsilon$ {\normalfont then}\\
\mbox{}\hspace{6ex} $x_{m+1}:=z$\\
\mbox{}\hspace{6ex} $\mathcal{D}:=(x_i)_{i=1}^{m+1}$, $\mathcal{I}:=\mathcal{I}\cup\{m+1\}$\\
\mbox{}\hspace{6ex} $m:=m+1$\\
\mbox{}\hspace{3ex} {\normalfont else}\\
\mbox{}\hspace{6ex} $x_{m+1}:=z$, $x_{m+2}:=x^*(z)$\\
\mbox{}\hspace{6ex} $\mathcal{D}:=(x_i)_{i=1}^{m+2}$, $\mathcal{E}:=\mathcal{E}\cup\{m+1\}$, $\mathcal{B}:=\mathcal{B}\cup\{m+2\}$\\
\mbox{}\hspace{6ex} $m:=m+2$.\\
\mbox{}\hspace{3ex} {\normalfont endif}\\
{\normalfont end}
\end{Algorithm}

The idea behind this algorithm is simple.
For any $z\in\tilde\Omega$, the fact that $\theta(z)=0$ implies $z\in\mathcal{R}(T,t_0,x_0)$.
If, on the other hand, $\theta(z)>0$, then $z\notin\mathcal{R}(T,t_0,x_0)$ and $x^*(z)\in\partial\mathcal{R}(T,t_0,x_0)$,
assuming no error has been made in the global optimisation routine.
This way, Algorithm \ref{construct:training:set} constructs a training set $\mathcal{D}:=(x_i)_{i=1}^m$
of points $(x_i)_{i=1}^m\subset\Omega$ with index set partitions $\mathcal{I}$, $\mathcal{E}$ and $\mathcal{B}$. % and labels $(y_i)_{i=1}^m\subset\{-1,1\}$ with
%\[y_i=\left\{\begin{array}{rl}1, & x_i\in\mathcal{R}(T,t_0,x_0)\\ -1, & \text{otherwise}\end{array}\right..\]
%$y_i=1$ if and only if $x_i\in\mathcal{R}(T,t_0,x_0)$.
The small parameter $\epsilon>0$ is introduced to compensate for numerical precision errors.
By construction, we have $|\mathcal{I}| + |\mathcal{E}| + |\mathcal{B}| = m$.

The support vector machine algorithm is designed to find a function from an RKHS (along with its sublevel set)
which best fits a labelled training set such as this, in the sense of minimising a suitable loss function.
However, the context here differs from the usual setting in which the SVM is applied (to a set of randomly
generated data potentially subject to noise) in two main ways.

Firstly, the training set is not just labelled according to whether a sample point belongs to the reachable set or not,
but also has the possible label of being on the boundary, as indicated by the three index sets $\mathcal{I}$, $\mathcal{E}$
and $\mathcal{B}$.

Secondly, the standard soft-margin SVM classifier allows for statistical errors in the labelling of data.
Here, only specific errors can occur: a point that is labelled as an interior or boundary point must belong to the reachable set
as the optimal control routine finds an admissible path to reach that point. So we do not want to allow a point
with index $i\in\mathcal{I}\cup\mathcal{B}$ to be on the exterior of our reachable set approximation. However, a point labelled as an exterior point
could well belong to the reachable if the optimal control routine failed to find the global minimum.

 We present the following adapted SVM algorithm in order to account for these differences:

\begin{Algorithm}[Adapted SVM] \label{AdaptedSVM}
First run Algorithm \ref{construct:training:set} to obtain the set $\mathcal{D} = (x_i)_{i=1}^m$ and index set partitions $\mathcal{I}$, $\mathcal{E}$ and $\mathcal{B}$. Fix regularisation parameters $C_1,C_2>0$, and let $K(\cdot,\cdot)$ be a Mercer %symmetric, positive-definite
kernel with corresponding finite-dimensional RKHS $\mathcal{H}_K$ on $X:=\{x_1,x_2,\dots,x_m\}$.
We search for a function $f=\sum_{i=1}^m a_i K_{x_i}$  in $\mathcal{H}_K$
by solving the following optimisation problem over the optimisation variables $(a, b, \xi, \eta) \in \mathbb{R}^m\times\mathbb{R} \times \mathbb{R}^{|\mathcal{E}|} \times \mathbb{R}^{|\mathcal{B}|}$:
\begin{eqnarray}
\operatorname{Minimise }_{a,b,\xi,\eta}  \quad \frac{1}{2}\|f\|_{\mathcal{H}_K}^2
+ C_1\sum_{i\in\mathcal{E}} \xi_i + C_2\sum_{i\in\mathcal{B}}\eta_i, &	\label{AdaptedSVM:cost}\\
\text{subject to }\qquad \sum_{k=1}^m a_k K(x_k,x_i) + b  \ge 1,& \qquad i\in\mathcal{I},	\label{AdaptedSVM:Iconstraints}\\
-\sum_{k=1}^m a_kK(x_k,x_i) - b \ge 1 - \xi_i,& \qquad i\in\mathcal{E},	\label{AdaptedSVM:Econstraints}\\
\sum_{k=1}^ma_kK(x_k,x_i) + b = \eta_i,& \qquad i\in\mathcal{B},	\label{AdaptedSVM:Bconstraints}\\
\xi_i\ge 0,& \qquad i\in\mathcal{E},	\label{AdaptedSVM:Eslack}\\
\eta_i \ge 0,& \qquad i\in\mathcal{B}.	\label{AdaptedSVM:Bslack}
\end{eqnarray}
The approximation of the reachable set $\mathcal{R}(T,t_0,x_0)$ is given by
\begin{equation} \label{eqn:SVMrepresentation}
\mathcal{R}_{SVM}(T,t_0,x_0):=\{x\in\mathbb{R}^d : f(x) + b \ge 0\}.
\end{equation}
\end{Algorithm}

The labelled training set generated by Algorithm \ref{construct:training:set}, which contains all available
knowledge about the reachable set, is incorporated in constraints \eqref{AdaptedSVM:Iconstraints}, \eqref{AdaptedSVM:Econstraints}
and \eqref{AdaptedSVM:Bconstraints}.
The constraint \eqref{AdaptedSVM:Iconstraints} ensures that the function value is at least $1$ on the points
that are labelled as interior points.
Note that there is no slack variable appearing in this constraint, according to our observation that points labelled
as interior points must lie within the reachable set.
In contrast, \eqref{AdaptedSVM:Econstraints} contains the non-negative slack variable $\xi_i$ (see also \eqref{AdaptedSVM:Eslack}),
which allows for the possibility of an error being made on a point labelled as an exterior point.
Where the slack variable is zero, the function value is less than or equal to $-1$ on exterior points.
The constraint \eqref{AdaptedSVM:Bconstraints} tries to place boundary points on the level set $\{x\in\mathbb{R}^d:f(x)+b=0\}$.
Here the non-negativity condition \eqref{AdaptedSVM:Eslack} follows front the fact that points labelled as boundary points
are the endpoints of orbits of \eqref{svrkm} and so cannot be on the exterior of the reachable set.

The first term of the cost function \eqref{AdaptedSVM:cost} controls the complexity of the function $f\in\mathcal{H}_K$
(and hence the sub level set) to avoid overfitting the training set $\mathcal{D}$. % through minimising the norm of the function.
This is contrasted with the following two terms, which control the penalty due to errors in classification.
This bias-variance trade-off is managed through the regularisation coefficients $C_1$ and $C_2$.
As these coefficients approach infinity, the function $f$ is allowed to become more and more complex, and the solution to the
optimisation problem approaches the hard-margin solution where no errors are permitted on the training set.

The optimisation problem \eqref{AdaptedSVM:cost}--\eqref{AdaptedSVM:Bslack} is a convex optimization problem, and in particular all the constraints are affine. In this case Slater's Theorem guarantees strong duality if the problem is feasible \cite{Boyd_04_1}. In the case of the Gaussian kernel, feasibility is guaranteed by the following theorem \cite{Micchelli_86_1}.

\begin{theorem}   \label{thm:Gaussianfullrank}
Let $x_1,\ldots,x_m\subset X$ be distinct points, and $\sigma >0$. The matrix $\mathbb{K}$ given by
\begin{equation*}
\mathbb{K}_{ij} = \exp\big({-\textstyle\frac{1}{\sigma}\|x_i-x_j\|^2}\big)
\end{equation*}
has full rank.
\end{theorem}

Therefore, Algorithm \ref{AdaptedSVM} can be recast into the dual problem using the KKT conditions.
This is the problem that is generally solved in practice.

We introduce the variables $y_i \in\{-1,1\}$, $i=1,\ldots,m$ by defining $y_i = 1$ for $i\in\mathcal{I}\cup\mathcal{B}$ and $y_i=-1$ for $i\in\mathcal{E}$.

\begin{Algorithm}[Dualised SVM]
Under the same conditions as in Algorithm \ref{AdaptedSVM}, solve the following minimisation problem over the variables $\alpha \in \mathbb{R}^m$:
\begin{eqnarray}
\operatorname{Minimise }_{\alpha}  \quad \frac{1}{2}\sum_{i,j=1}^m y_i y_j \alpha_i  \alpha_j K(x_i,x_j) - \sum_{i\in \mathcal{I}\cup\mathcal{E}}\alpha_i &	\label{AdaptedSVMdual:cost}\\
\text{subject to }\qquad \sum_{i=1}^m y_i\alpha_i = 0,	\label{AdaptedSVMdual:offsetconstraints}\\
\alpha_i\ge 0,\qquad i\in\mathcal{I}	\label{AdaptedSVMdual:Iconstraints}\\
0\le\alpha_i \le C_1, \qquad i\in\mathcal{E}	\label{AdaptedSVMdual:Econstraints}\\
\alpha_i \ge -C_2,\qquad i\in\mathcal{B}   \label{AdaptedSVMdual:Bconstraints}
\end{eqnarray}
\end{Algorithm}

The solution to the problem \eqref{AdaptedSVMdual:cost}--\eqref{AdaptedSVMdual:Bconstraints} provides the function $f=\sum_{i=1}^m a_i K_{x_i}$, where the $a_i$ are given by $a_i = y_i \alpha_i$ for $i\in\{1,\ldots, m\}$. The points $x_i$ for which the corresponding constraint \eqref{AdaptedSVMdual:Iconstraints}, \eqref{AdaptedSVMdual:Econstraints} or \eqref{AdaptedSVMdual:Bconstraints} are strictly satisfied are called the \textit{support vectors} in the literature. For the support vectors the corresponding constraints \eqref{AdaptedSVM:Iconstraints}, \eqref{AdaptedSVM:Econstraints} and \eqref{AdaptedSVM:Bconstraints} are satisfied as equalities, and in addition $\xi_i$ or $\eta_i$ is equal to zero. The offset $b$ can therefore be computed from \eqref{AdaptedSVM:Iconstraints}, \eqref{AdaptedSVM:Econstraints} and \eqref{AdaptedSVM:Bconstraints} for the support vectors.

Accordingly, points $x_i$ ($i\in\mathcal{E}$) for which $\alpha_i = C_1$ and points $x_i$ ($i\in\mathcal{B}$) for which $\alpha_i = C_2$  are the so-called \textit{error vectors}. These are the points for which $\xi_i$ and $\eta_i$ may be nonzero, and for which the reachable set approximation \eqref{eqn:SVMrepresentation} may misclassify. A boundary point $x_i$ for which $\eta_i >0$ will still be classified as being in the reachable set, but will not be on the boundary of the set approximation. However an exterior point $x_i$ for which $\xi_i>1$ will be misclassified by \eqref{eqn:SVMrepresentation}. If $0<\xi_i<1$ then $x_i$ will be still be on the exterior of $\mathcal{R}_{SVM}$ but will be inside the `margin' $\{x\in\mathbb{R}^d : |f(x)+b| <1\}$ and so it is still called an error vector.

Finally, points $x_i$ $(i\in\mathcal{I}\cup\mathcal{E})$ for which $\alpha_i=0$ are  \textit{ignored vectors}. It is not hard to see that these points have no influence on the solution to the above optimisation problem, and could as well have been left out of the data set. In addition, the property of being an ignored vector is robust with respect to perturbation of the support and error vectors. Note that the set of boundary points $x_i$ ($i\in\mathcal{B}$) by definition does not contain any ignored vectors, since the property that $\alpha_i=0$  is not robust with respect to such a perturbation due to \eqref{AdaptedSVM:Bconstraints}. Roughly speaking, points that are far away from the boundary of the reachable set (and are correctly classified) will be ignored vectors. However it is not practically possible to tell in advance which data points will be ignored vectors, or even if ignored vectors will remain ignored with the addition of new points.

%The expressions \textit{support vectors}, \textit{error vectors} and \textit{ignored vectors} are standard terminology
%in the machine learning community, cf. e.g. \cite{Cauwenberghs_01_1,Cortes_95_1}.

\subsection{Incremental updates} \label{sec:incremental}

It is possible to increase the accuracy of the SVM approximation
step by step, until a desired precision is reached.
In that case, the optimisation problem \eqref{AdaptedSVM:cost}--\eqref{AdaptedSVM:Bslack}
(or \eqref{AdaptedSVMdual:cost}--\eqref{AdaptedSVMdual:Bconstraints}) needs to be solved after  each addition of a batch of new points. This optimisation problem runs over all points in the training set, so as this set becomes larger, this may become costly.

Fortunately it is possible to solve the SVM optimisation problem by means of incremental updates \cite{Cauwenberghs_01_1}.
This procedure consists of deriving equations to keep the KKT conditions satisfied, as a new dual variable $\alpha_i$ is incremented
from zero.
The procedure ends when a new point becomes either a support vector or error vector. For details, we refer to \cite{Cauwenberghs_01_1}. Here, we outline the procedure for our adapted version of the SVM algorithm.

The cost function in the dual formulation of the optimisation problem \eqref{AdaptedSVMdual:cost} may trivially be rewritten in the more convenient form
\begin{equation}
 \frac{1}{2}\sum_{i,j=1}^m y_i y_j \alpha_i  \alpha_j K(x_i,x_j) - \sum_{i\in \mathcal{I}\cup\mathcal{E}}\alpha_i + b\sum_{i=1}^m y_i\alpha_i 	\label{AdaptedSVMdual:costinc}
\end{equation}
retaining the constraints \eqref{AdaptedSVMdual:offsetconstraints}--\eqref{AdaptedSVMdual:Bconstraints}, with the offset $b$ re-introduced as a Lagrange multiplier. The necessary and sufficient KKT conditions for this problem may be written as follows:
\begin{eqnarray}
i\in\mathcal{I}: & g_i:= f(x_i) + b -1 &
\begin{cases}
\ge 0 \quad \textrm{when } \alpha_i =0\\
=0 \quad \textrm{when } \alpha_i >0
\end{cases}		\label{eqn:KKTI}\\
i\in\mathcal{E}: & g_i:= -f(x_i) -b -1 &
\begin{cases}
\ge 0 \quad \textrm{when } \alpha_i =0\\
=0 \quad \textrm{when } 0 < \alpha_i < C_1\\
\le 0 \quad \textrm{when } \alpha_i = C_1
\end{cases}		\label{eqn:KKTE}\\
i\in\mathcal{B}: & g_i:= f(x_i) + b &
\begin{cases}
\ge 0 \quad \textrm{when } \alpha_i =-C_2\\
=0 \quad \textrm{when } \alpha_i > -C_2
\end{cases}		\label{eqn:KKTB}
\end{eqnarray}
\begin{equation}
\sum_{i=1}^m y_i\alpha_i  = 0				\label{eqn:KKTb}
\end{equation}

The conditions \eqref{eqn:KKTI}--\eqref{eqn:KKTB} are satisfied with equality for the support vectors. Given a new labelled point $(x_c,y_c)$ with dual variable $\alpha_c$ initially set to zero, we need to ensure that these equality conditions (as well as \eqref{eqn:KKTb}) continue to be satisfied for the support vectors as we increment $\alpha_c$ from zero.  Following \cite{Cauwenberghs_01_1}, we define the \textit{coefficient sensitivities} $\beta_i$ by
\begin{equation}		\label{eqn:betadef}
\left[
\begin{matrix}
\beta_0\\
\beta_{s_1}\\
\vdots\\
\beta_{s_{\mathcal{N}(S)}}
\end{matrix}
\right] = - \left[
\begin{matrix}
0 & y_{s_1} & \cdots & y_{s_{\mathcal{N}(S)}}\\
y_{s_1} & Q_{s_1s_1} & \cdots & Q_{s_1s_{\mathcal{N}(S)}}\\
\vdots & \vdots & \ddots & \vdots\\
y_{s_{\mathcal{N}(S)}} & Q_{s_{\mathcal{N}(S)}s_1} & \cdots &Q_{s_{\mathcal{N}(S)}s_{\mathcal{N}(S)}}
\end{matrix}
\right]^{-1}
\left[
\begin{matrix}
y_c\\ Q_{s_1c}\\
\vdots \\ Q_{s_{\mathcal{N}(S)}c}
\end{matrix}
 \right]
\end{equation}
where $Q_{ij}=y_iy_jK(x_i,x_j)$ and $\{s_1,\ldots,s_{\mathcal{N}(S)}\}$ is the index set corresponding to the support vectors ($\mathcal{N}(S)$ is the number of support vectors). We define $\beta_i=0$ for indices $i$ corresponding to ignored and error vectors. Then, the KKT conditions \eqref{eqn:KKTI}--\eqref{eqn:KKTb} will continue to be satisfied as $\alpha_c$ is incremented from zero provided the existing dual coefficients are also incremented according to
\begin{eqnarray}
\Delta b & = & \beta_0 \Delta\alpha_c	\label{eqn:deltab}\\
\Delta \alpha_i & = & \beta_i \Delta\alpha_c,	\label{eqn:deltaalpha}
\end{eqnarray}
where $\Delta\alpha_c$ is a small increment in $\alpha_c$. The \textit{margin sensitivities} $\gamma_i$ are likewise defined by
\begin{equation}
\gamma_i = Q_{ic} + \sum_{j=s_1}^{s_{\mathcal{N}(S)}}Q_{ij}\beta_{j} + y_i\beta_0
\end{equation}
and give the variation of the margins $g_i$ in \eqref{eqn:KKTI}--\eqref{eqn:KKTB}:
\begin{equation}
\Delta g_i = \gamma_i\Delta\alpha_c	\label{eqn:deltag}
\end{equation}
Equations \eqref{eqn:betadef}--\eqref{eqn:deltaalpha} ensure that $\gamma_i = 0$ for support vectors.

Now, for each new point $(x_c,y_c)$ the corresponding $g_i$ is first computed. If the new point automatically satisfies the KKT conditions then it is an ignored vector and $\alpha_c$ is left at zero. Otherwise we use \eqref{eqn:deltab}, \eqref{eqn:deltaalpha} and \eqref{eqn:deltag} to compute the largest possible increment of $\alpha_c$ so that \eqref{eqn:KKTI}--\eqref{eqn:KKTb} continue to be satisfied, at which point either $x_c$ becomes a support or error vector, or else another point in the data set migrates between the sets of support, error or ignored vectors. Then the coefficient and margin sensitivities must be recomputed and the procedure continues.

\subsection{Ball checking} \label{sec:ball:check}

The procedure outlined in the previous can also naturally be reversed in order to remove a point from the training set.
In practical implementation we have included the ball checking routine as described in the end of section \ref{section:dfog}
in which case it is sometimes necessary to remove an exterior point from the data set.
Note however that in these cases the SVM algorithm allows to keep the corresponding boundary point in the data set.

In practice, where the optimisation routine fails to find the global minimum, it is usually the case that the computed optimal
trajectory still terminates at a boundary point of the reachable set.
However, note that the SVM algorithm does allow for points in the index set $\mathcal{B}$ to actually be interior points.
In this case, some information on the reachable set is still retained in the case of an error due to the global optimisation
routine of the DFOG method.

\section{Examples}

We illustrate the qualities of our method by applying it to two examples from the literature.
In each example we compare its performance in the reachable set representation $\mathcal{R}_{SVM}(T,t_0,x_0)$
with that of the DFOG method and its reachable set approximation
$\mathcal{R}_{DFOG}(T,t_0,x_0)$ as given in \eqref{eqn:DFOGrepresentation}.

\subsection{A bilinear control system}

The following example is taken from \cite{Baier_13_1, Hajek_08_1} as a model system that exhibits convexity of the reachable set
for small times, but nonconvexity for larger times.
Both the DFOG and SVM methods work for either case.
%Our numerical approximations are for the reachable set at time $T=1$, where the reachable set is nonconvex.

We consider the two-dimensional control system
\begin{align}
\dot{x}_1 & =  \pi x_2\,,	\label{eqn:bilinear1}\\
\dot{x}_2 & =  -\pi u(t)x_1\,,	\label{eqn:bilinear2}\\
x_1(0) & =  -1\,,	\label{eqn:bilinear3}\\
x_2(0) & =  0\,,	\label{eqn:bilinear4}\\
u(t)  \in U & =  [0,1]\,.	\label{eqn:bilinear5}
\end{align}

\begin{figure}[ht]
	       \begin{center}
	        \includegraphics[width=0.7\textwidth, natwidth=311,natheight=159]{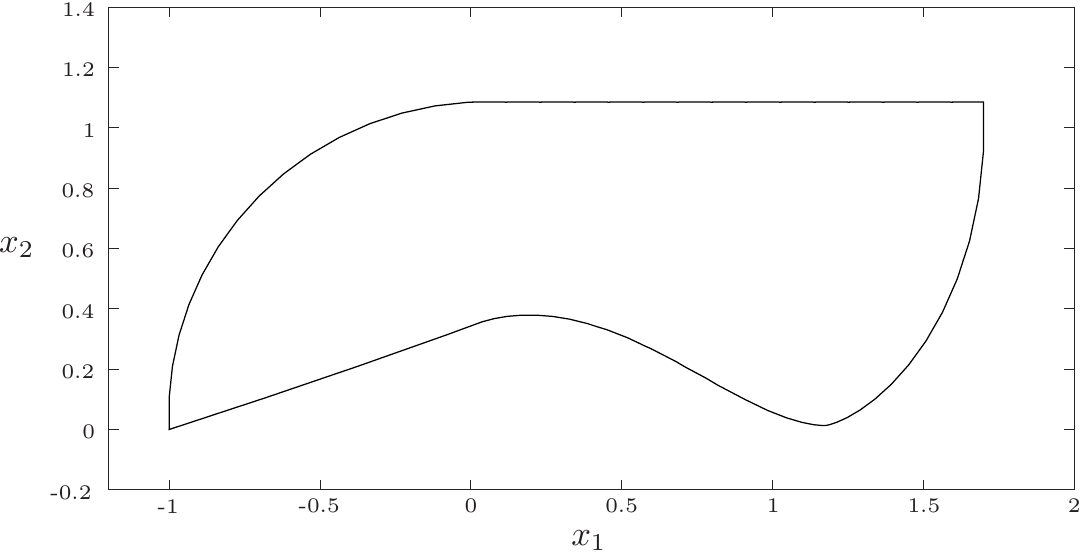}
	        \end{center}
\caption{Reachable set $\mathcal{R}(1,0,(-1,0))$ for the bilinear control system.}
\label{fig:bilinearDFOG}
\end{figure}

We are interested in approximating the reachable set $\mathcal{R}(1,0,{x}_0)$ for ${x}_0:=(x_1(0),x_2(0))$.
The reachable set is shown in Figure \ref{fig:bilinearDFOG}.
In this computation the time interval $[0,1]$ has been discretised with $N=30$ steps. We note that the error due to time discretisation is the same for both the DFOG and modified SVM methods. This is because both methods use the same time discretisation in the constraints for the optimal control problem from Algorithm \ref{DFOGmethod}. The difference between the two methods is the spatial representation of the reachable set.
In order to compare the methods we leave the time discretisation at $N=30$ and vary the spatial grid size $\rho$.

\begin{figure}
\centering
\begin{subfigure}[b]{0.45\textwidth}
	       \raisebox{0.1cm}{\includegraphics[width=\textwidth, natwidth=311,natheight=158]{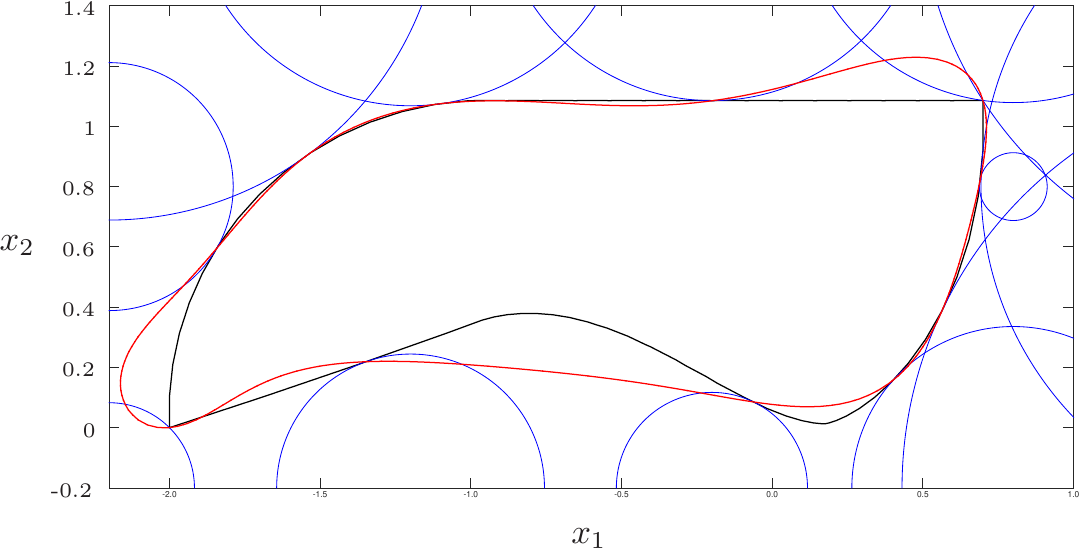}}
	        %\vspace{-0.3cm}
                \caption{$\rho = 1$}
                \label{fig:BLrho=1}
        \end{subfigure}
        \hspace{0.2cm}
             \begin{subfigure}[b]{0.45\textwidth}
                 \raisebox{0.1cm}{\includegraphics[width=\textwidth, natwidth=311,natheight=158]{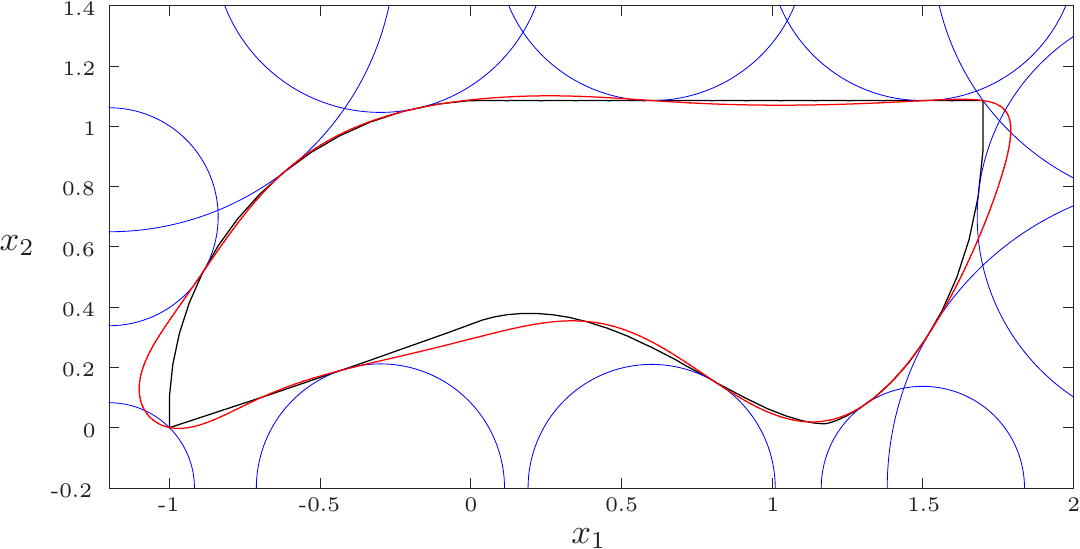}}
                \caption{$\rho = 0.9$}
                \label{fig:BLrho=0_9}
        \end{subfigure}
        %%%%%%%%%%%%
        \vspace{0.4cm}
                \hfill
        \begin{subfigure}[b]{0.45\textwidth}
	      \raisebox{0.1cm}{\includegraphics[width=\textwidth, natwidth=311,natheight=158]{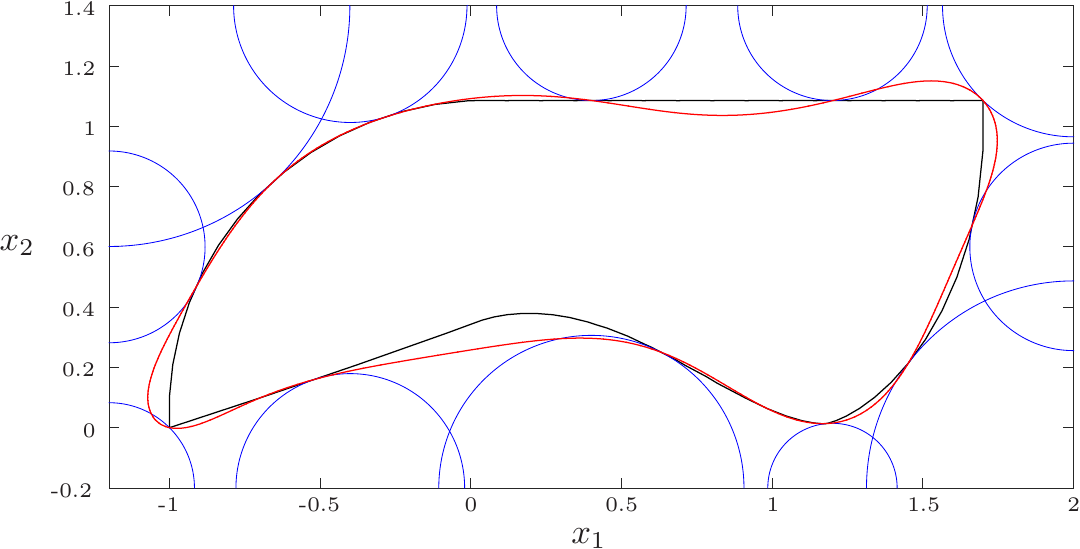}}
                \caption{$\rho = 0.8$}
                \label{fig:BLrho=0_8}
        \end{subfigure}
         \hspace{0.2cm}
             \begin{subfigure}[b]{0.45\textwidth}
             \raisebox{0.1cm}{\includegraphics[width=\textwidth, natwidth=311,natheight=158]{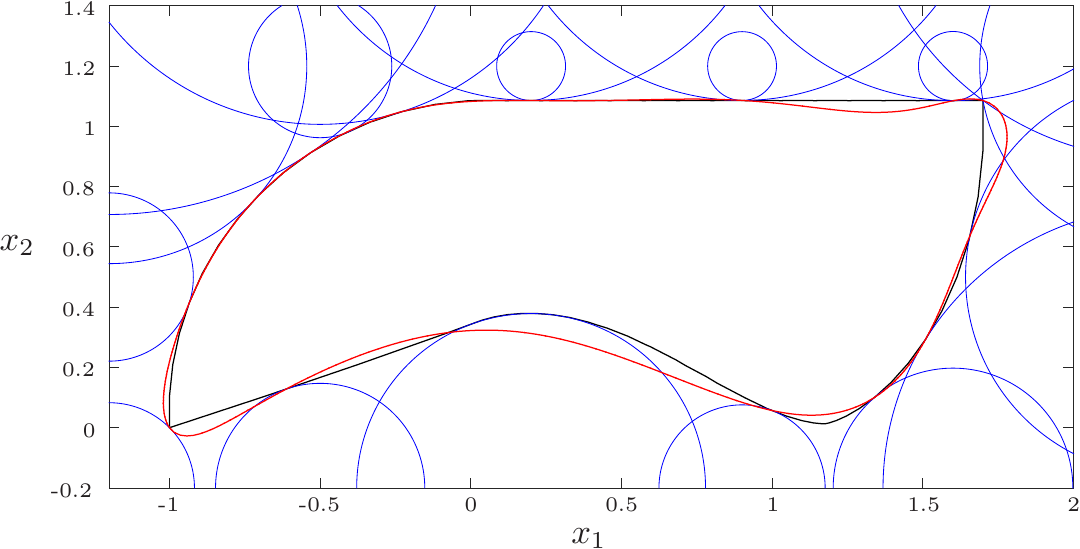}}
                \caption{$\rho = 0.7$}
                \label{fig:BLrho=0_7}
        \end{subfigure}
                %%%%%%%%%%%%%%
        \vspace{0.4cm}
                \hfill
        \begin{subfigure}[b]{0.45\textwidth}
	     \raisebox{0.1cm}{\includegraphics[width=\textwidth, natwidth=311,natheight=158]{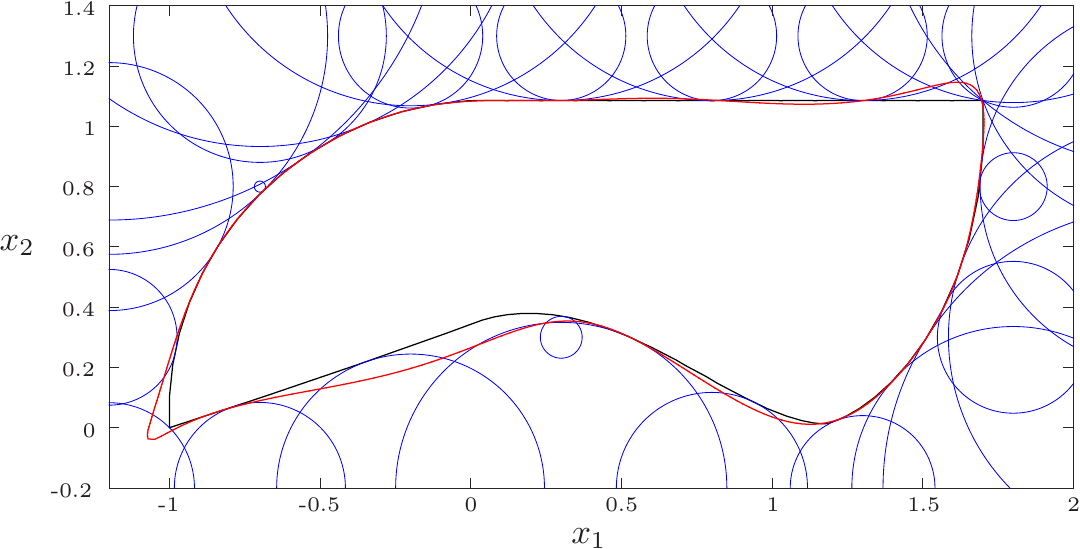}}
                \caption{$\rho = 0.5$}
                \label{fig:BLrho=0_5}
        \end{subfigure}
             \hspace{0.2cm}
             \begin{subfigure}[b]{0.45\textwidth}
              \raisebox{0.1cm}{\includegraphics[width=\textwidth, natwidth=311,natheight=158]{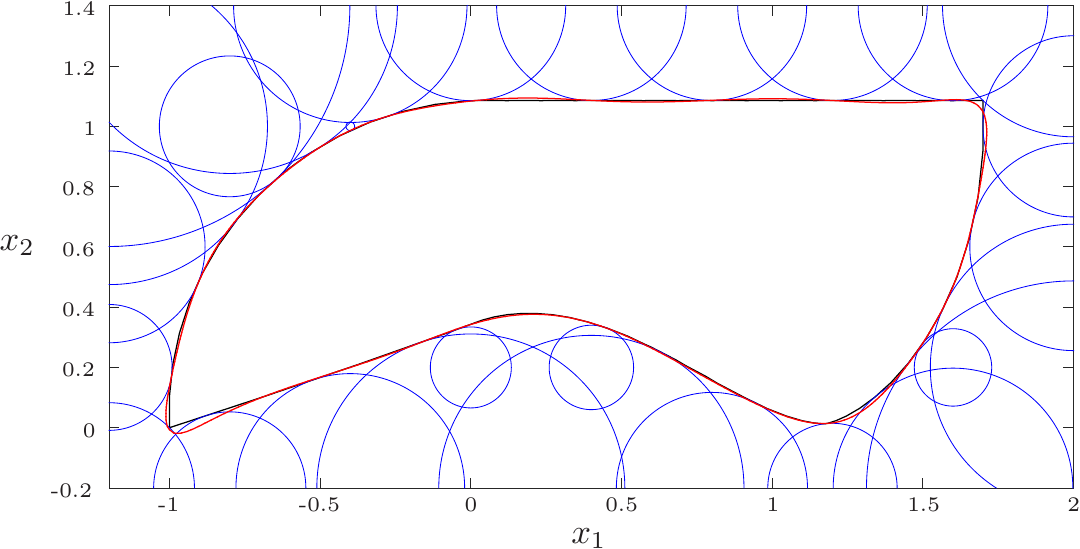}}
                \caption{$\rho = 0.4$}
                \label{fig:BLrho=0_4}
        \end{subfigure}
            %%%%%%%%%%%%%%
        \vspace{0.4cm}
                \hfill
        \begin{subfigure}[b]{0.45\textwidth}
	   \raisebox{0.1cm}{\includegraphics[width=\textwidth, natwidth=311,natheight=158]{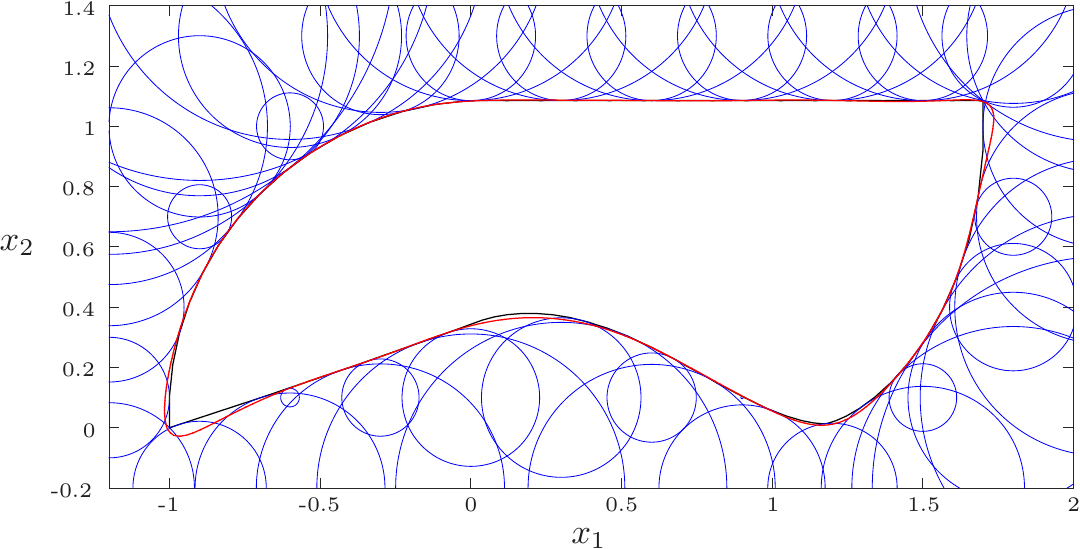}}
                \caption{$\rho = 0.3$}
                \label{fig:BLrho=0_3}
        \end{subfigure}
             \hspace{0.2cm}
             \begin{subfigure}[b]{0.45\textwidth}
            \raisebox{0.1cm}{\includegraphics[width=\textwidth, natwidth=311,natheight=158]{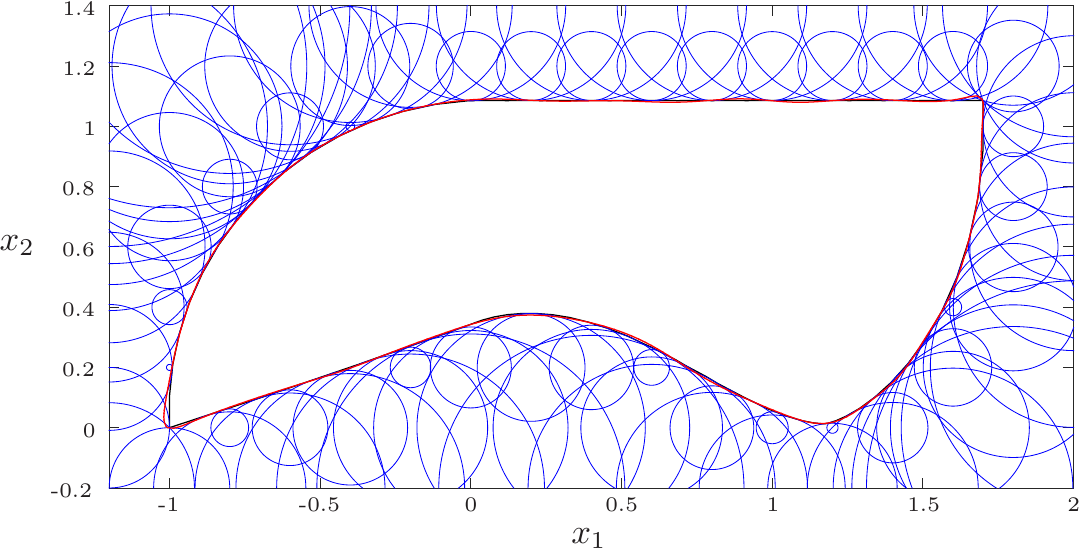}}
                \caption{$\rho = 0.2$}
                \label{fig:BLrho=0_2}
        \end{subfigure}
\caption{Reachable set approximation for the bilinear control system, $T=1$, $N=30$. The exact reachable set is shown in blue, with the SVM approximations shown in red and the DFOG approximations shown as green circles. Approximations are shown for varying grid sizes $\rho$.}
\label{fig:bilinearDFOG_SVM}
\end{figure}

Figure \ref{fig:bilinearDFOG_SVM} shows successive approximations of both the DFOG and SVM methods for the reachable set $\mathcal{R}_h(1,0,{x}_0)$, where $h=\frac{1}{N}$ and $N=30$ is fixed. For both algorithms, the set of grid points $\tilde{\Omega}$ was defined as a restriction of $\rho\mathbb{Z}^2$, and the approximations are made  for varying spatial discretizations $\rho$ independently. In this example there are very few errors made by the global optimisation routine. The Hausdorff distances calculated between the true reachable set and the numerical approximations from both the DFOG and SVM methods are shown in Figure \ref{fig:BLCS_HD}.

\begin{minipage}[c]{\textwidth}
\begin{minipage}[c]{0.38\textwidth}
\begin{center}
\begin{tabular}{ccc}
\hline\hline
$\rho$ & DFOG & SVM\\[0.5ex]
\hline
1.0 & 0.3794 & 0.1889\\
0.9 & 0.3738 & 0.1177\\
0.8 & 0.3780 & 0.1034\\
0.7 & 0.2373 & 0.0792\\
0.6 & 0.2165 & 0.0919\\
0.5 & 0.1542 & 0.0803\\
0.4 & 0.1113 & 0.0296\\
0.3 & 0.0604 & 0.0412\\
0.2 & 0.0222 & 0.0190\\[1ex]
\hline
\end{tabular}
\end{center}
\end{minipage}
\begin{minipage}[c]{0.6\textwidth}
%\centering
\begin{center}
\includegraphics[width=0.75\textwidth, natwidth=323,natheight=250]{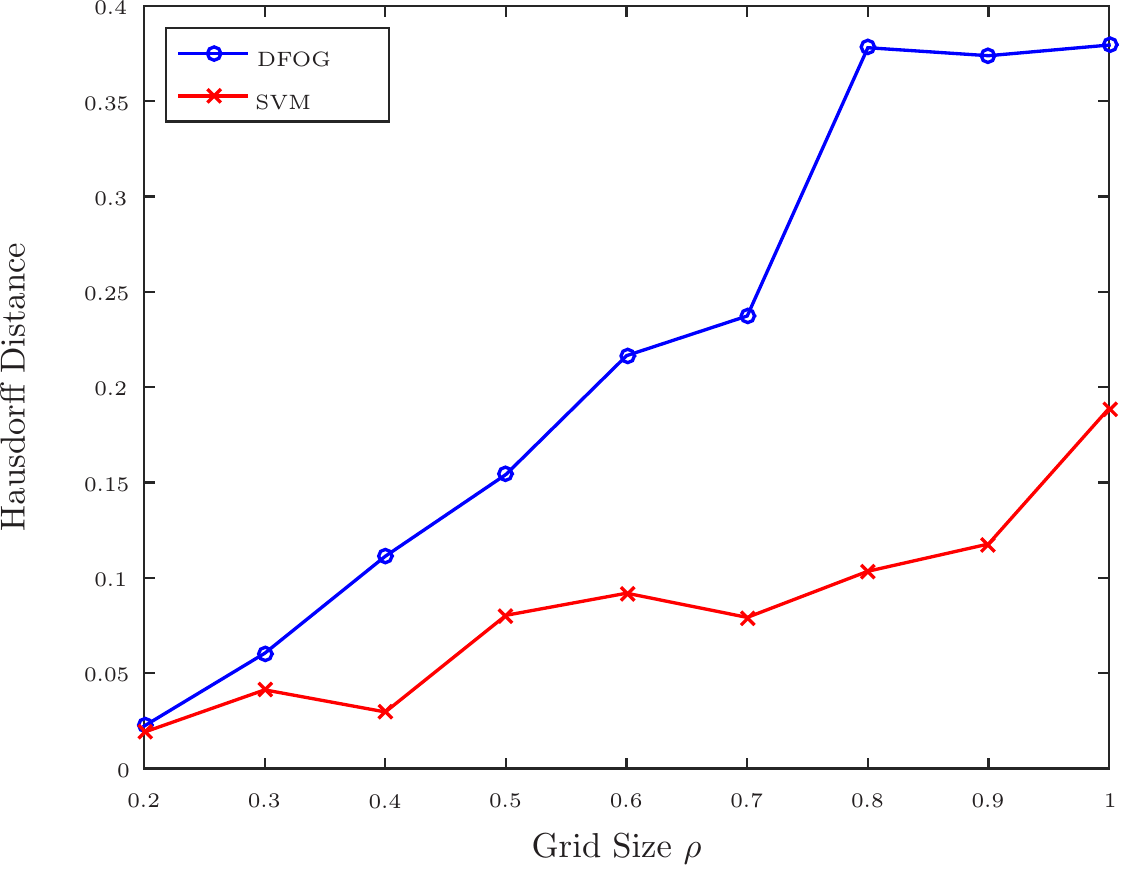}
\end{center}
\end{minipage}
\captionof{figure}{Table and Figure showing Hausdorff distance values calculated between the true reachable set and SVM and DFOG approximations for varying grid sizes $\rho$.}
\label{fig:BLCS_HD}
\end{minipage}

\subsection{A nonlinear control system}

The following example was presented in \cite{Rieger_Unpub_2} as an example of a reachable set that may change its topology for different times $T$.
\begin{align}
\dot{x}_1 & =  x_1(1-|x_1|) - x_1x_2 + u_1\,,\\
\dot{x}_2 & =  x_1^4 - \tfrac{1}{2} + u_2\,,\\
x_1(0) & = 0\,,\\
x_2(0) & =  0\,,\\
(u_1(t), u_2(t)) \in U & = [-\tfrac{1}{5},\tfrac{1}{5}]\times[-\tfrac{1}{5},\tfrac{1}{5}]\,.
\end{align}

\begin{figure}[h]
	        \begin{center}
\includegraphics[width=0.4\textwidth, natwidth=304,natheight=402]{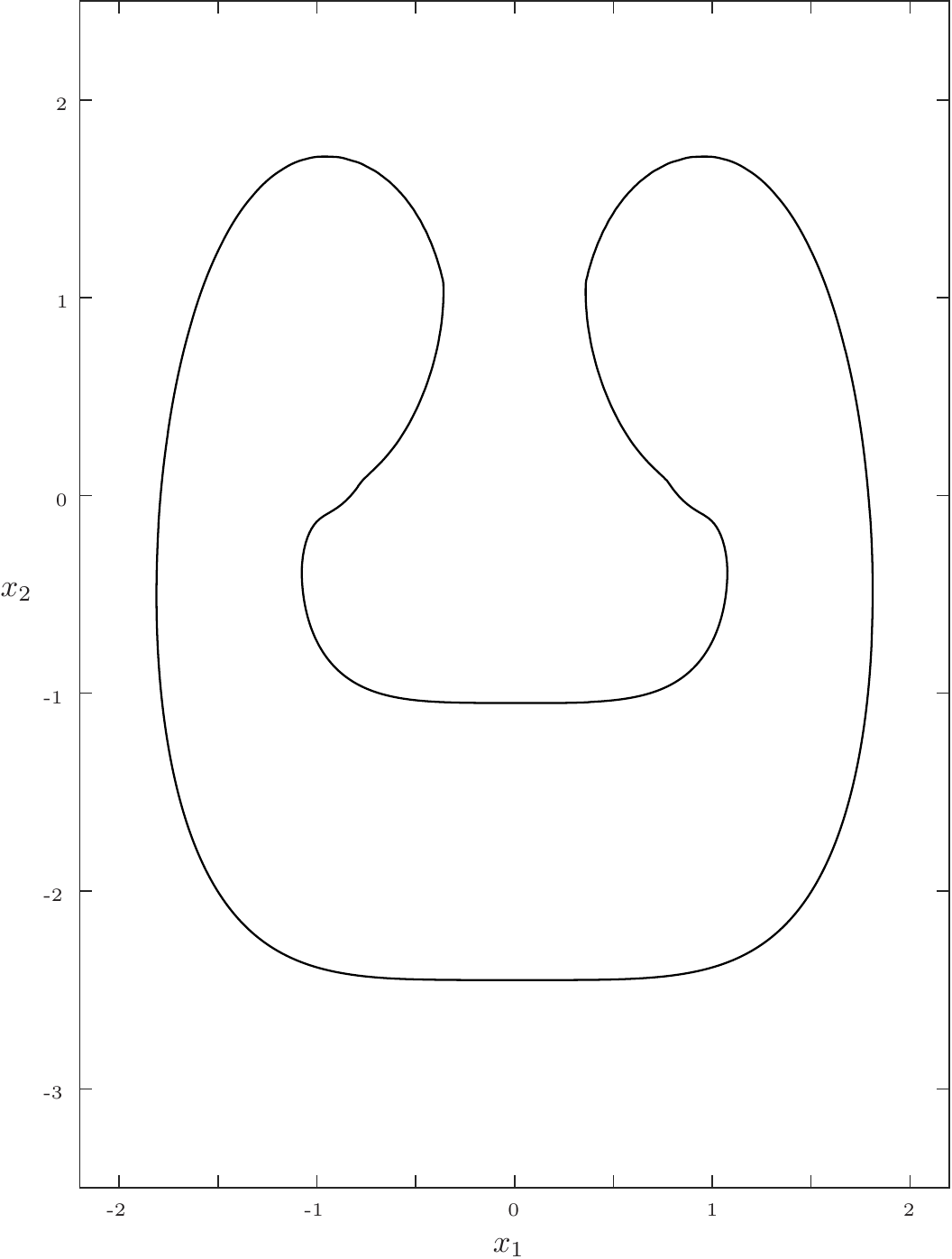}
	        \end{center}
\caption{Reachable set $\mathcal{R}(3.5,0,(0,0))$ for the nonlinear control system.}
\label{fig:topDFOG}
\end{figure}

In this example we will approximate the reachable set $\mathcal{R}(3.5,0,{x}_0)$, where ${x}_0 := (x_1(0),x_2(0))$. This reachable set is shown in Figure \ref{fig:topDFOG}.

\begin{figure}
\centering
\begin{subfigure}[b]{0.28\textwidth}
	      \raisebox{-0cm}{\includegraphics[width=\textwidth, natwidth=304,natheight=402]{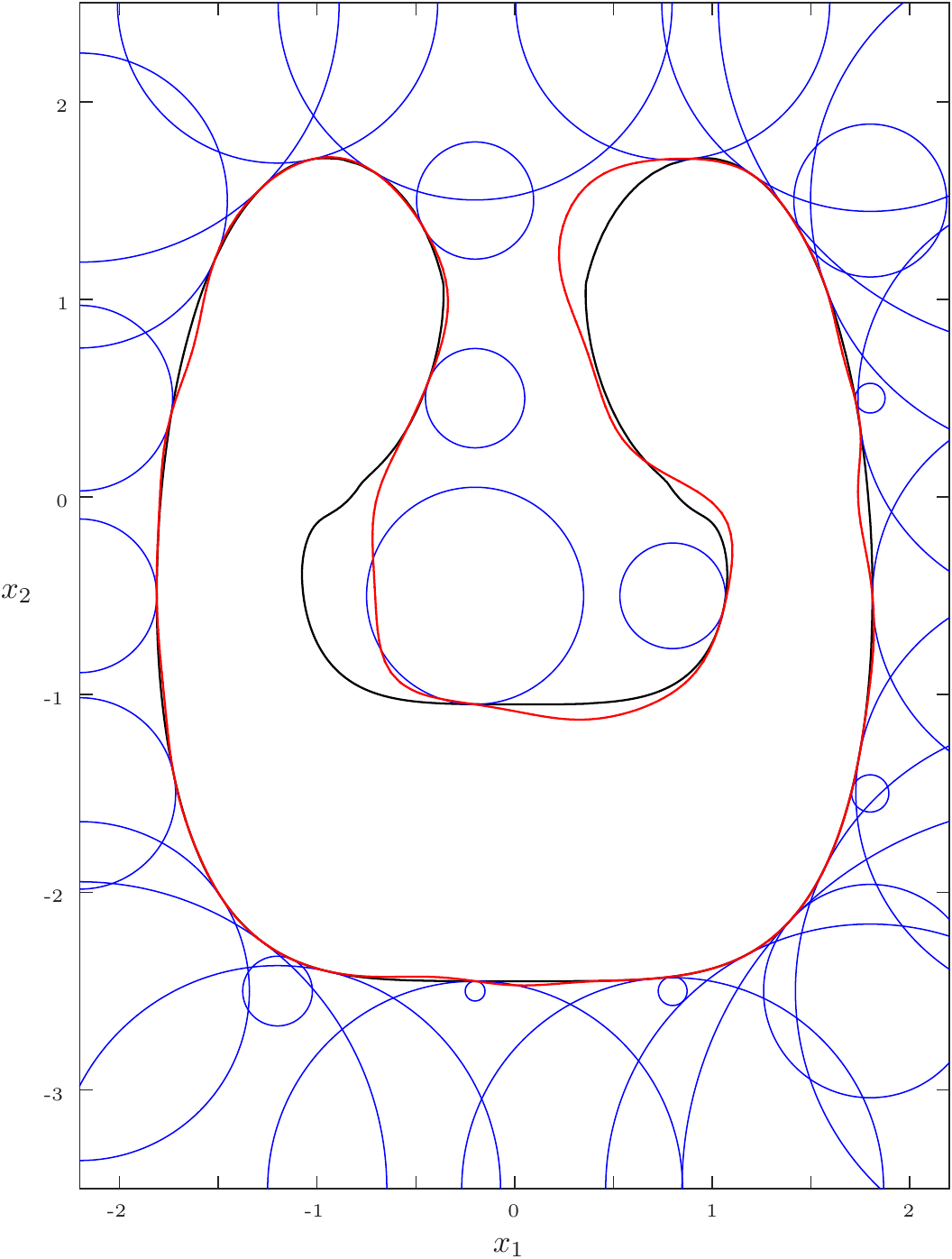}}
	        \vspace{-0.65cm}
                \caption{$\rho = 1$}
                \label{fig:NLrho=1}
        \end{subfigure}
        \hspace{0.4cm}
             \begin{subfigure}[b]{0.28\textwidth}
            \raisebox{-0cm}{\includegraphics[width=\textwidth, natwidth=304,natheight=402]{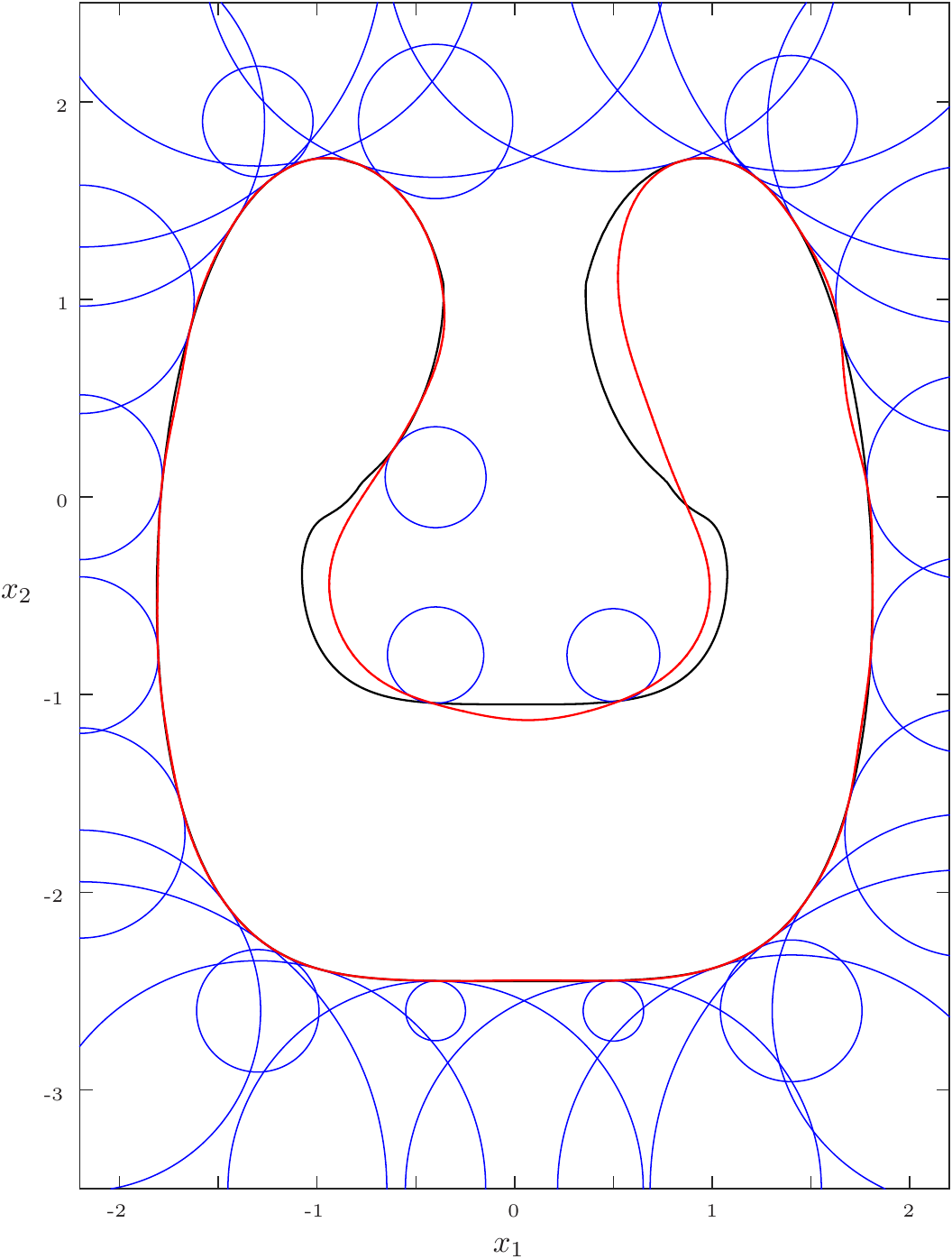}}
                 \vspace{-0.65cm}
                \caption{$\rho = 0.9$}
                \label{fig:NLrho=0_9}
        \end{subfigure}
        \hspace{0.4cm}
        \begin{subfigure}[b]{0.28\textwidth}
          \raisebox{0cm}{\includegraphics[width=\textwidth, natwidth=304,natheight=402]{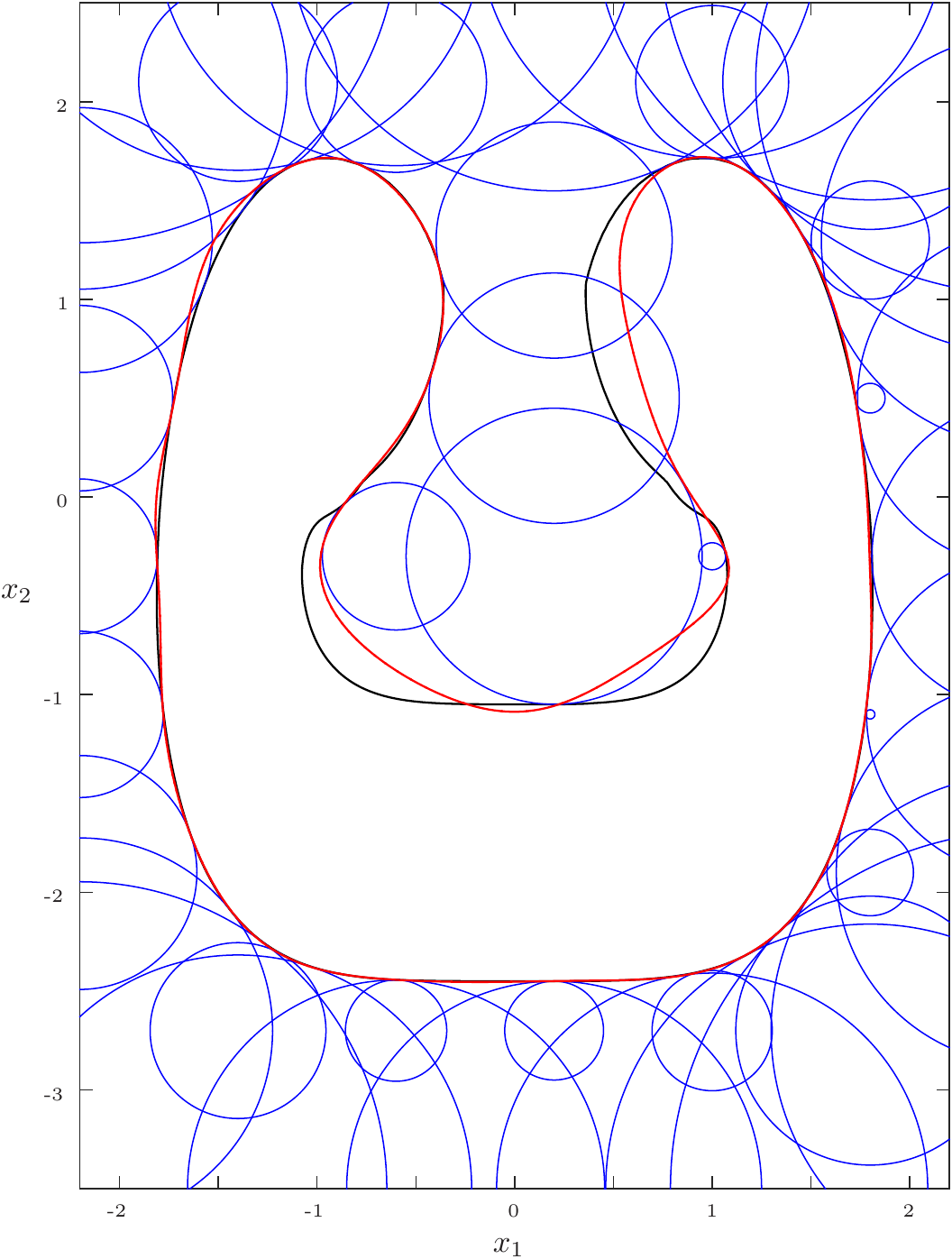}}
                 \vspace{-0.65cm}
                \caption{$\rho = 0.8$}
                \label{fig:NLrho=0_8}
        \end{subfigure}
        %%%%%%%%%%%%
        \vspace{-0cm}
                \hfill
        \begin{subfigure}[b]{0.28\textwidth}
         \raisebox{0cm}{\includegraphics[width=\textwidth, natwidth=304,natheight=402]{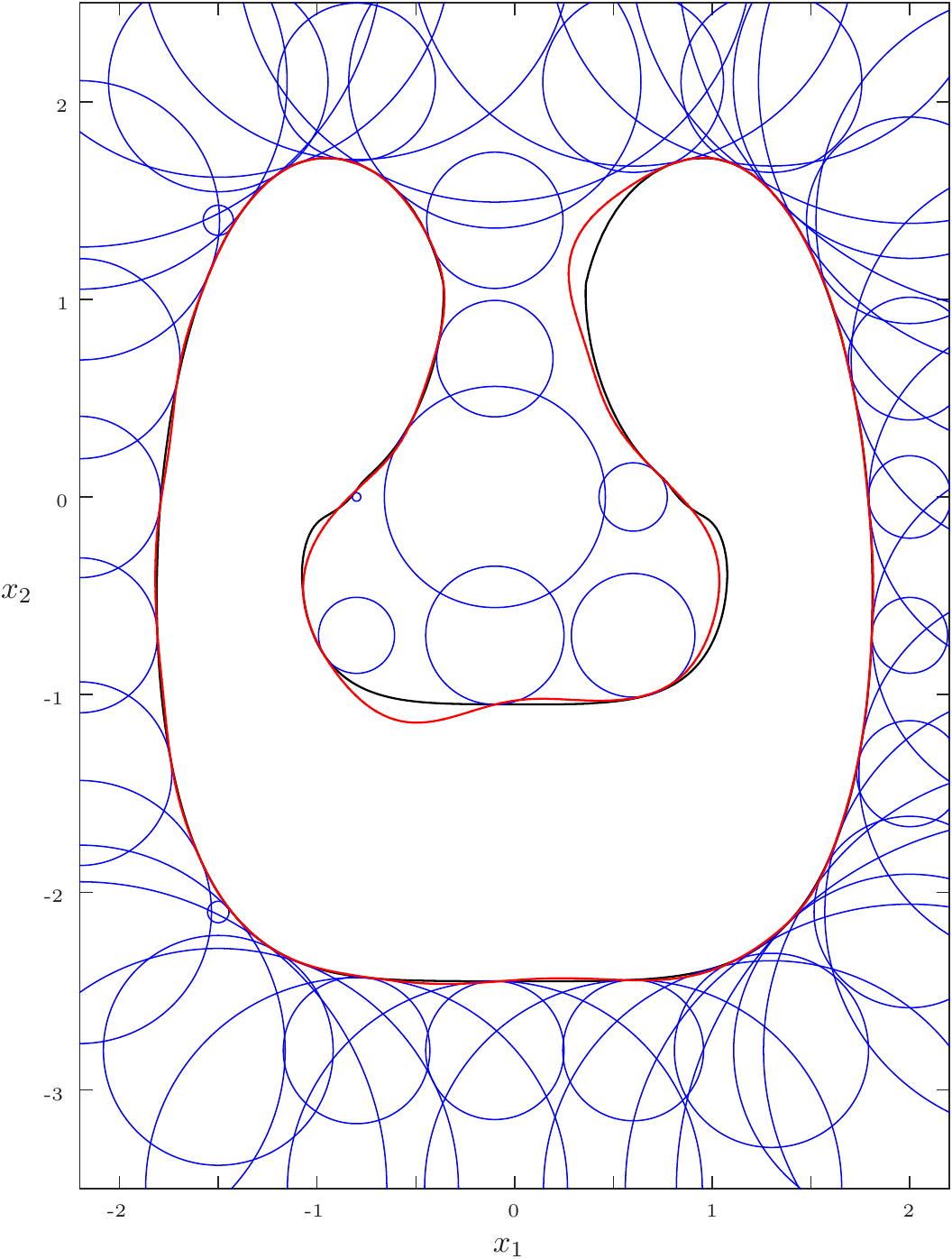}}
                 \vspace{-0.65cm}
                \caption{$\rho = 0.7$}
                \label{fig:NLrho=0_7}
        \end{subfigure}
        \hspace{0.4cm}
             \begin{subfigure}[b]{0.28\textwidth}
             \raisebox{0cm}{\includegraphics[width=\textwidth, natwidth=304,natheight=402]{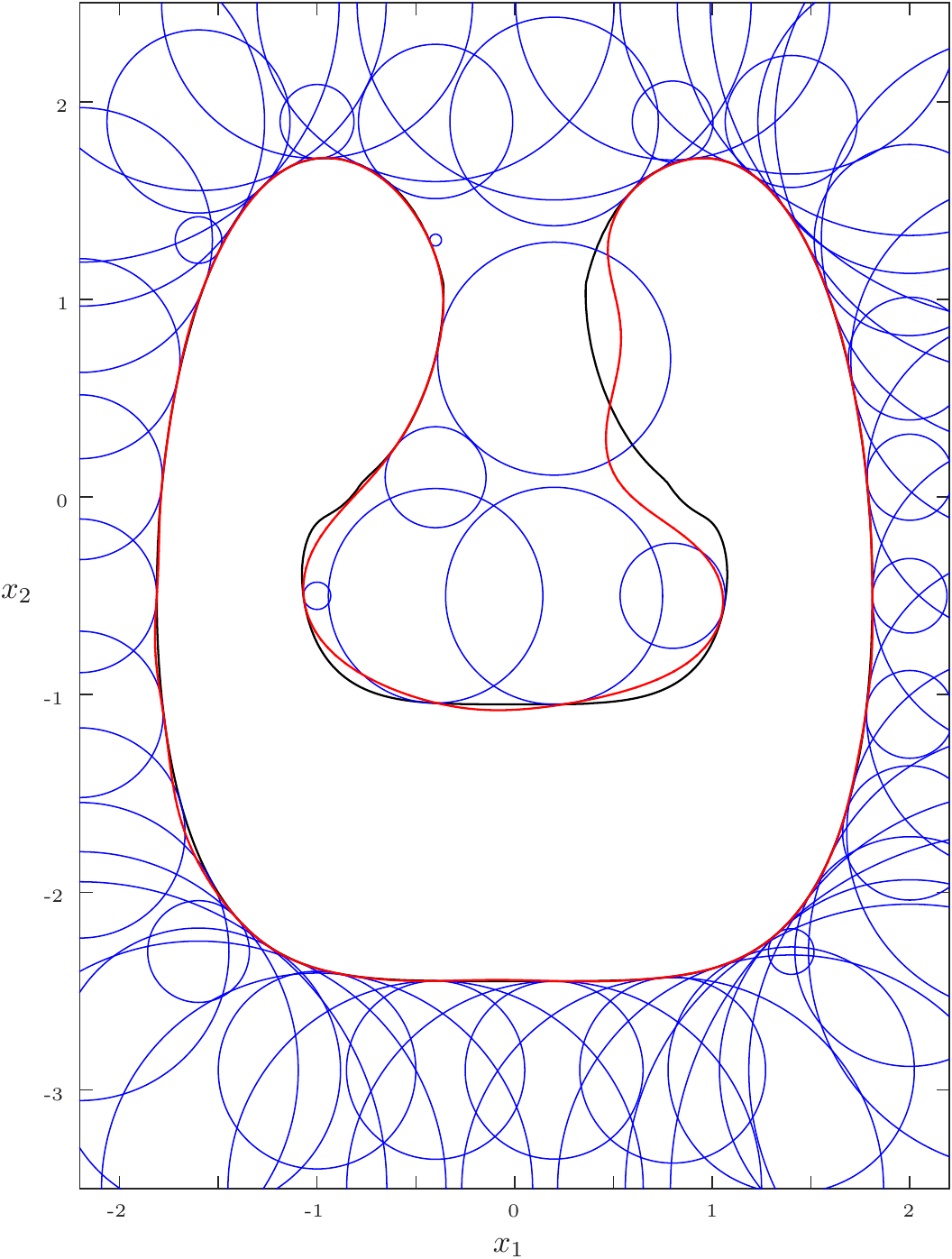}}
                 \vspace{-0.65cm}
                \caption{$\rho = 0.6$}
                \label{fig:NLrho=0_6}
        \end{subfigure}
        \hspace{0.4cm}
        \begin{subfigure}[b]{0.28\textwidth}
         \raisebox{0cm}{\includegraphics[width=\textwidth, natwidth=304,natheight=402]{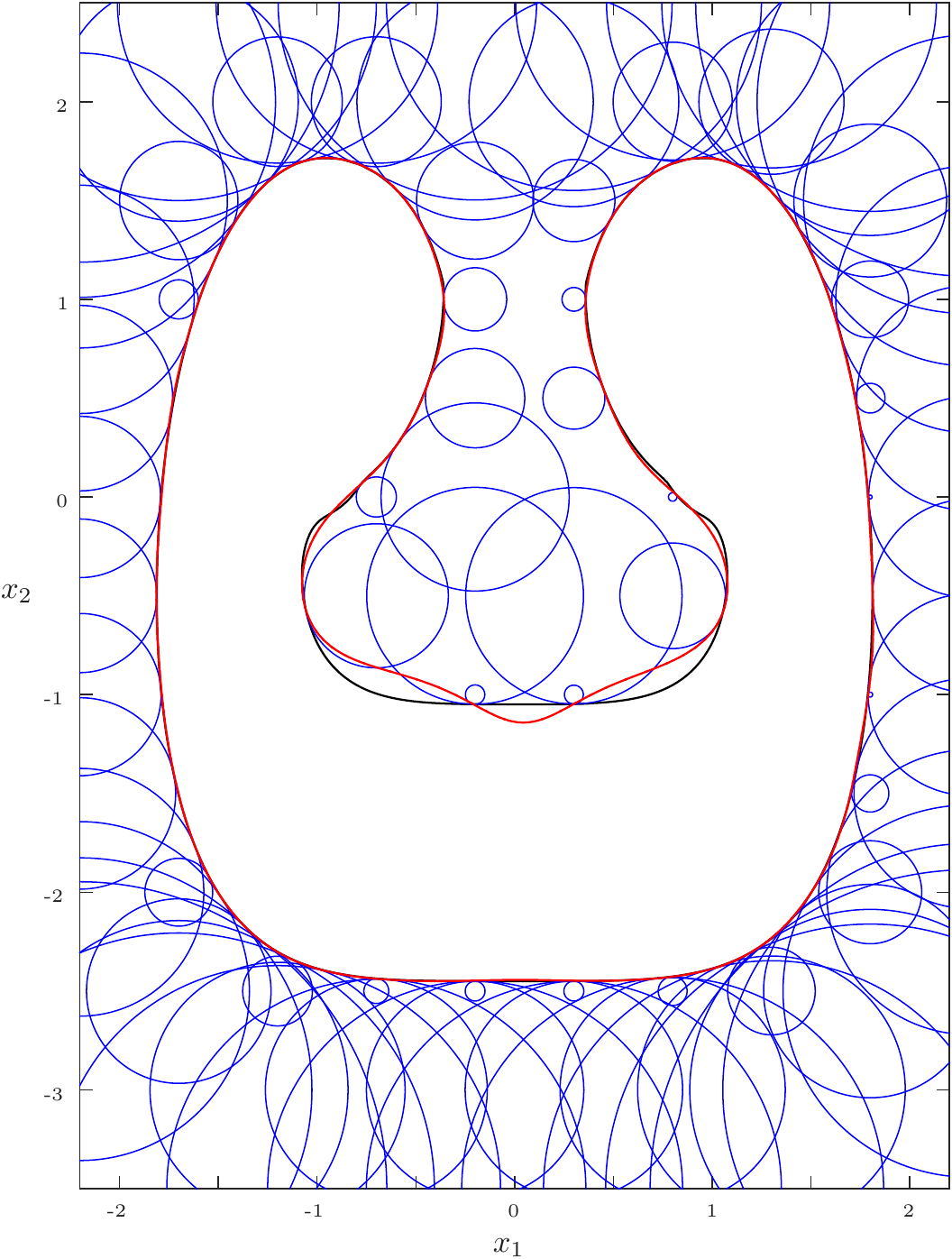}}
                 \vspace{-0.65cm}
                \caption{$\rho = 0.5$}
                \label{fig:NLrho=0_5}
         \end{subfigure}
                %%%%%%%%%%%%%%
        \vspace{-0cm}
                \hfill
           \begin{subfigure}[b]{0.28\textwidth}
           \raisebox{0cm}{\includegraphics[width=\textwidth, natwidth=304,natheight=402]{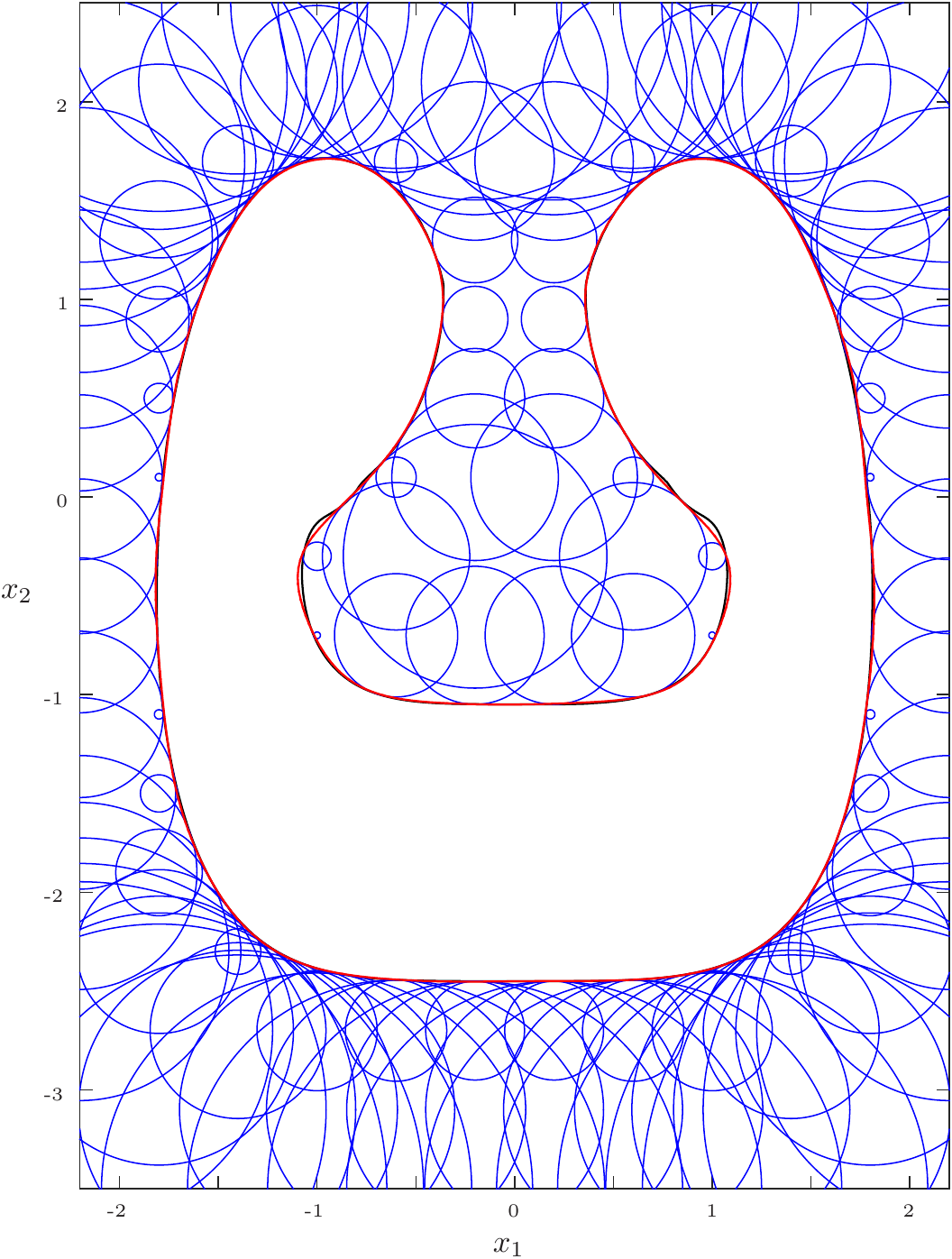}}
                 \vspace{-0.65cm}
                \caption{$\rho = 0.4$}
                \label{fig:NLrho=0_4}
        \end{subfigure}
        \hspace{0.4cm}
             \begin{subfigure}[b]{0.28\textwidth}
            \raisebox{0cm}{\includegraphics[width=\textwidth, natwidth=304,natheight=402]{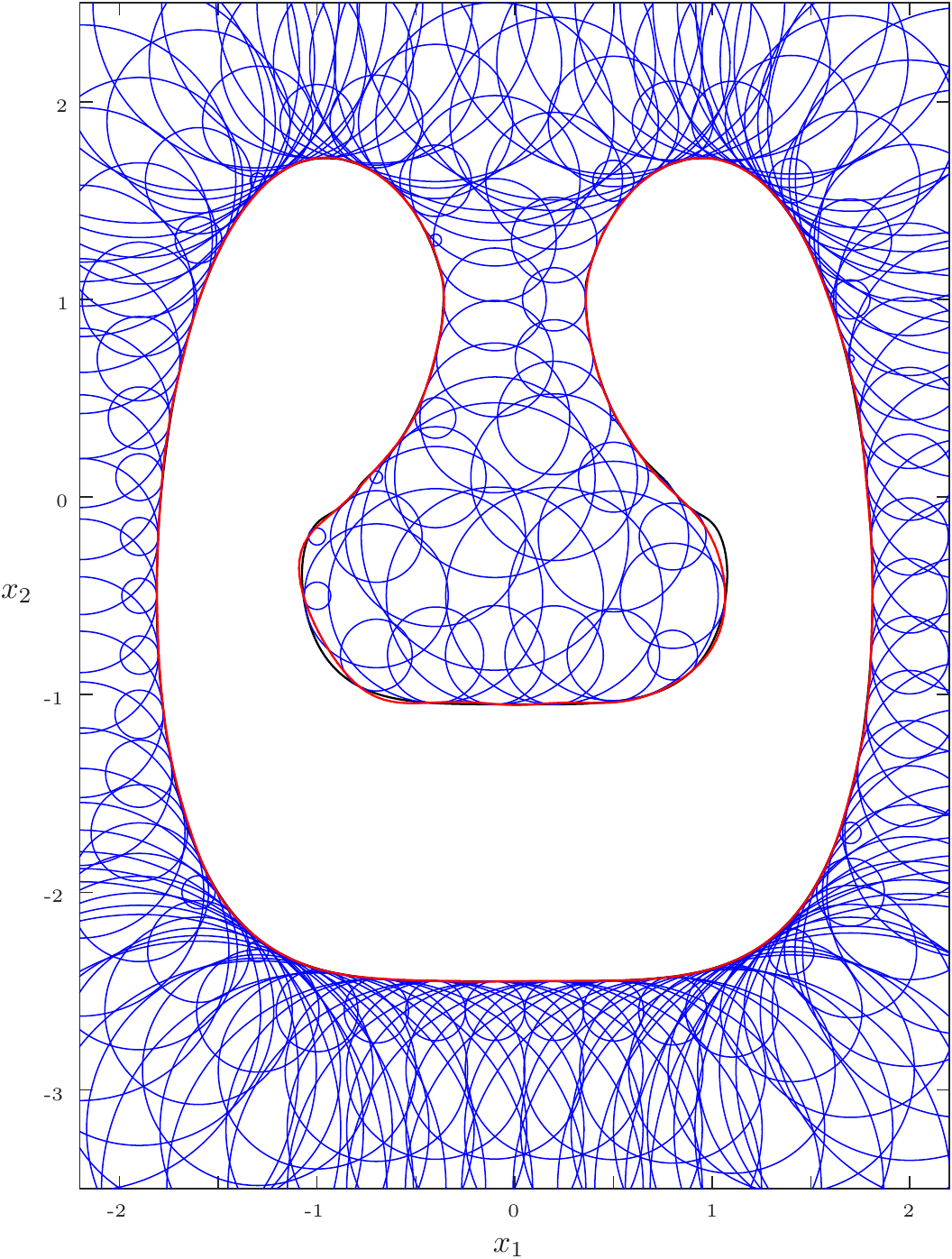}}
                \vspace{-0.65cm}
                \caption{$\rho = 0.3$}
                \label{fig:NLrho=0_3}
        \end{subfigure}
        \hspace{0.4cm}
        \begin{subfigure}[b]{0.28\textwidth}
       \raisebox{0cm}{\includegraphics[width=\textwidth, natwidth=304,natheight=402]{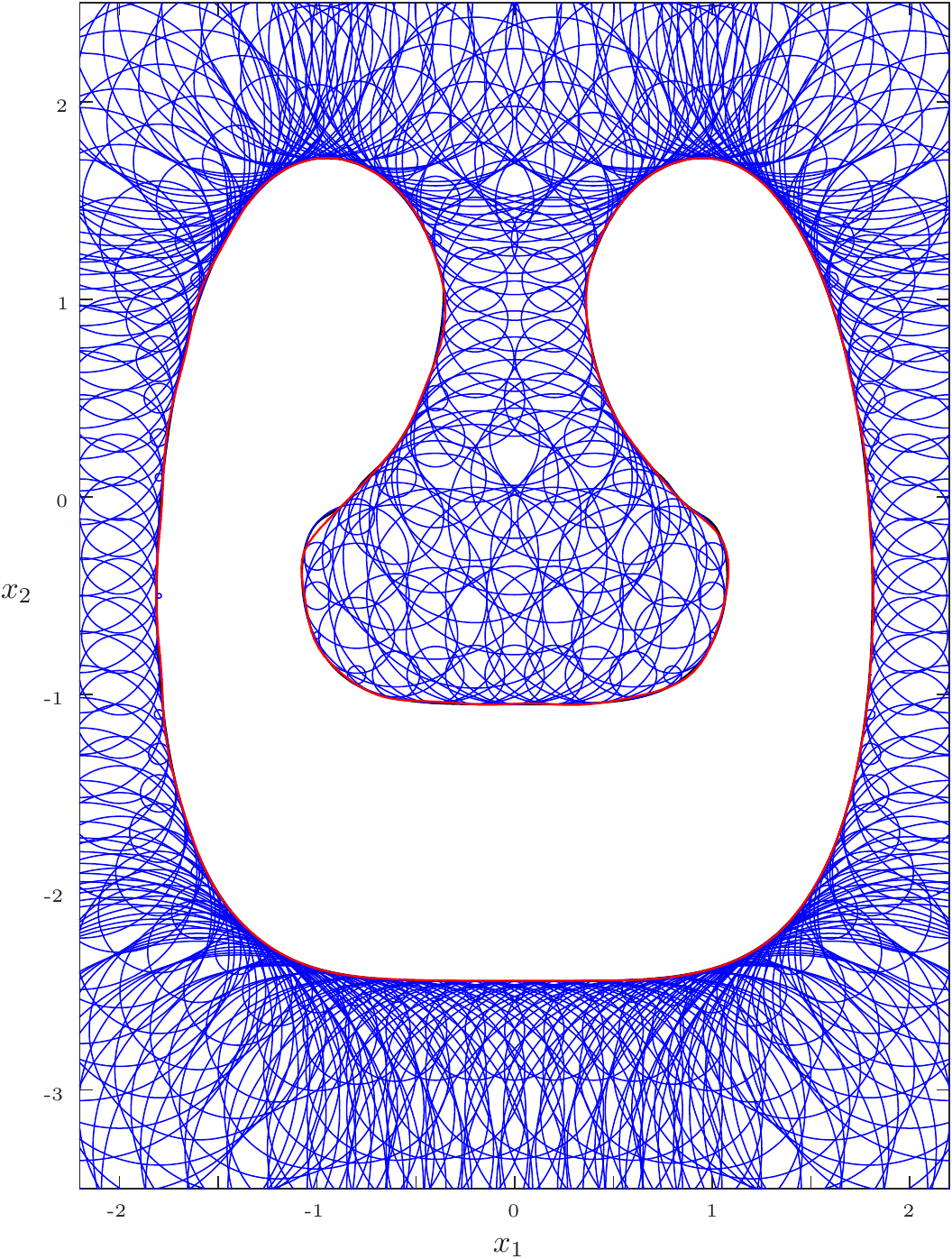}}
                \vspace{-0.65cm}
                \caption{$\rho = 0.2$}
                \label{fig:NLrho=0_2}
        \end{subfigure}
\caption{Reachable set approximation for the nonlinear control system, $T=3.5$, $N=100$. The exact reachable set is shown in blue, with the SVM approximations shown in red and the DFOG approximations shown as green circles. Approximations are shown for varying grid sizes $\rho$.}
\label{fig:RiegerDFOG_SVM}
\end{figure}

Figure \ref{fig:RiegerDFOG_SVM} shows successive approximations of  the DFOG and SVM methods for the reachable set $\mathcal{R}_h(3.5,0,{x}_0)$, for $N=50$. As before, the grid points are defined as a restriction of the grid $\rho\mathbb{Z}^2$. The approximations are shown for varying spatial discretizations $\rho$, and each figure is produced by an independent run of the algorithms.
Figure \ref{fig:NLCS_HD} provides the Hausdorff distances calculated between the true reachable set and the approximations made by both the SVM and DFOG methods.

This example contains more errors made by the global optimisation routine than the previous example, due to the highly non-convex topology of the reachable set. We can see that the SVM algorithm appears to converge faster to a good approximation of the reachable set. Again, the SVM algorithm is somewhat robust to these global optimisation errors. Note that information is still added to the SVM algorithm even in the case of a global optimisation error, since the computed optimal point will still be in the reachable set, and this point is still added to the algorithm.

\begin{minipage}[c]{\textwidth}
\begin{minipage}[c]{0.38\textwidth}
\begin{center}
\begin{tabular}{ccc}
\hline\hline
$\rho$ & DFOG & SVM\\[0.5ex]
\hline
1.0 & 0.6452 & 0.3629\\
0.9 & 0.7868 & 0.1862\\
0.8 & 0.3910 & 0.1994\\
0.7 & 0.6550 & 0.1128\\
0.6 & 0.6378 & 0.1755\\
0.5 & 0.5578 & 0.1256\\
0.4 & 0.4039 & 0.0321\\
0.3 & 0.1089 & 0.0443\\
0.2 & 0.0335 & 0.0288\\[1ex]
\hline
\end{tabular}
\end{center}
\end{minipage}
\begin{minipage}[c]{0.6\textwidth}
%\centering
\begin{center}
\includegraphics[width=0.75\textwidth, natwidth=320,natheight=250]{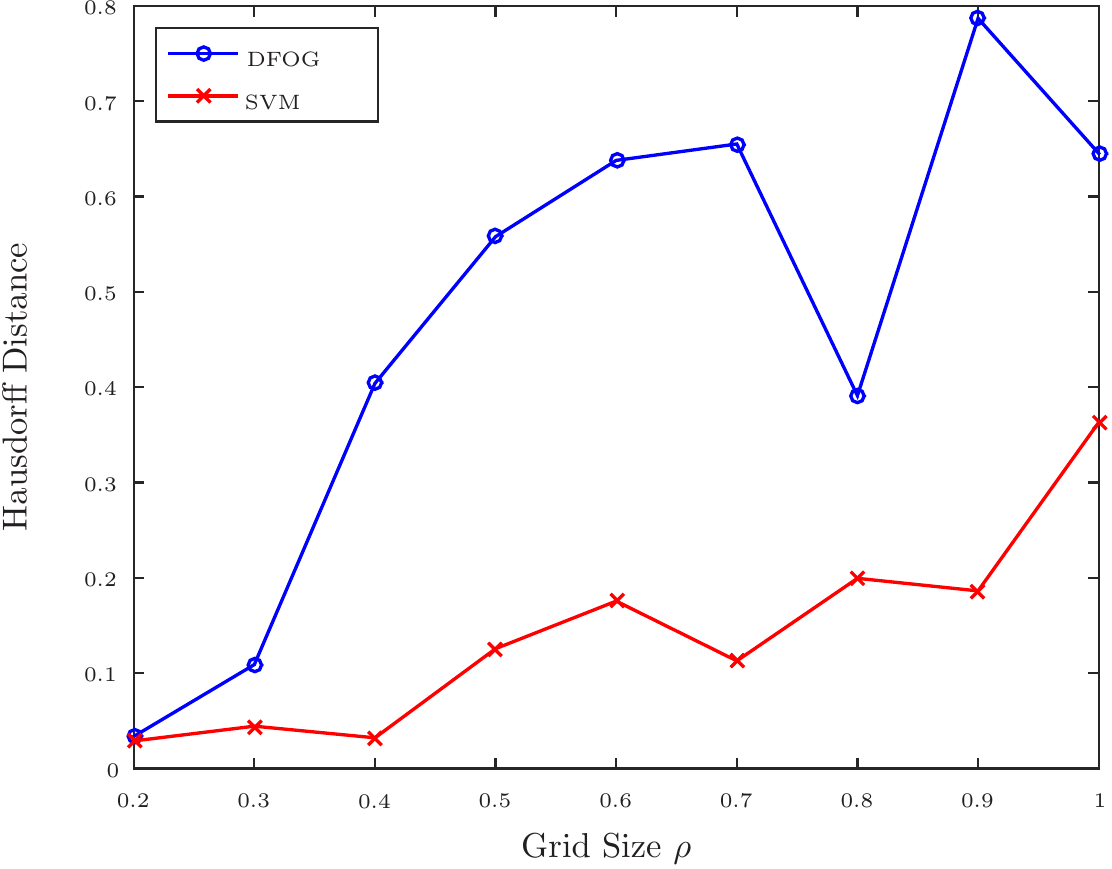}
\end{center}
\end{minipage}
\captionof{figure}{Table and Figure showing Hausdorff distance values calculated between the true reachable set and SVM and DFOG approximations for varying grid sizes $\rho$.}
\label{fig:NLCS_HD}
\end{minipage}

\section{Conclusions}
The modified Support Vector Machine algorithm provides an alternative representation of the reachable set, based on the results gained from a set of global optimisation problems provided by the DFOG algorithm. This new approach has the advantage that it is robust to a small number global optimisation errors, and appears to benefit from faster convergence for particular examples.
Several specialised algorithms exist for efficiently solving the standard SVM optimisation problem \cite{Joachims_98_1,Platt_98_1}, which could be adapted to the modified SVM algorithm we have presented here. The  global optimal control problems are also particularly expensive  when the dimension of the control variable is large, or when a fine time discretisation is used. Therefore for many real-world problems only few optimal control problems can be solved in practice. In these cases where relatively few data from the optimal control routine are
available, the SVM algorithm performs significantly better. In addition, the sublevel set representation of the SVM approach is more handy for many applications than the DFOG representation.

As is always the case with algorithms of this type, there are several parameters in the algorithm that need to be tuned for optimal performance. The tolerance $\epsilon$ as described in Section~\ref{sec:SVM} is important to distinguish interior and exterior points, and affects the approximation for both the DFOG and SVM algorithms. Within the SVM algorithm, the parameters $C_1$ and $C_2$ control the regularity of the solution as described earlier. Also for radial basis functions such as the Gaussian kernel used in our examples, the scaling parameter $\sigma$ is an additional parameter, related to the regularisation parameters $C_1$ and $C_2$. These parameters can be chosen using standard validation techniques such as hold-out testing (see \cite{Cauwenberghs_01_1,Vapnik_95_1}), but the precise effect of these parameters on the regularity of the solution (and how the parameters relate to each other) is not yet well understood.

A final important problem is that of choosing the best points on which to run the global optimisation routine in order to improve the current approximation. In our problem setting, we are in the fortunate position of being able to choose any point to run the algorithm on at each step. This is in contrast to many applications of the Support Vector Machine, where the data is randomly generated from an unknown underlying distribution. In our example applications we have run the algorithm on a regular grid, however it is clear that this is not the optimal strategy. The question of how to choose the best point is likely to be related to problem of understanding the effect of the parameters in the algorithm, and again is a worthwhile subject of future work. The framework provided in this paper to incrementally update the SVM algorithm is also a precursor to such a strategy.

A further benefit of our proposed methodology is that it may also in principle be used to compute invariant sets for random dynamical systems \cite{Homburg_06_1,Lamb_Unpub_1}, as well as invariant sets for control systems \cite{Colonius_00_1}.

\paragraph{Acknowledgements.}
The first and the third author were supported by the UK Engineering and Physical Sciences Research Council (EPSRC).

\providecommand{\bysame}{\leavevmode\hbox to3em{\hrulefill}\thinspace}
\providecommand{\MR}{\relax\ifhmode\unskip\space\fi MR }
\providecommand{\MRhref}[2]{%
  \href{http://www.ams.org/mathscinet-getitem?mr=#1}{#2}
}
\providecommand{\href}[2]{#2}


\begin{thebibliography}{CYLW98}

\bibitem[Aro50]{Aronszajan_50_1}
N.~Aronszajn, \emph{Theory of reproducing kernels}, Transactions of the
  American Mathematical Society \textbf{68} (1950), 337--404.

\bibitem[BBCG07]{BaiBusChaGer07}
R.~Baier, C.~B\"uskens, I.A.~Chahma, and M.~Gerdts, \emph{Approximation of reachable sets by direct solution methods for optimal control problems}, Optim. Methods Softw. \textbf{22} No. 3 (2007), 433--452.

\bibitem[BG09]{Baier_09_1}
R.~Baier and M.~Gerdts, \emph{A computational method for non-convex reachable
  sets using optimal control}, Proceedings of the European Control Conference
  (ECC) 2009, Budapest (Hungary), August 23--26, 2009 (Budapest), EUCA, 2009,
  pp.~97--102.

\bibitem[BGX13]{Baier_13_1}
R.~Baier, M.~Gerdts, and I.~Xausa, \emph{Approximation of reachable sets using
  optimal control algorithms}, Numerical Algebra, Control and Optimization
  \textbf{3} (2013), no.~3, 519--548.

\bibitem[BR07]{Beyn_07_1}
W.-J. Beyn and J.~Rieger, \emph{Numerical fixed grid methods for differential
  inclusions}, Computing \textbf{81} (2007), no.~1, 91--106.

\bibitem[BR10]{Beyn_10_1}
\bysame,
\emph{The implicit Euler scheme for one-sided {L}ipschitz differential inclusions},
Disc.~Cont.~Dyn.~Sys.~B \textbf{14} (2010), no.~2, 409--428.

\bibitem[BV04]{Boyd_04_1}
S.~Boyd and L.~Vandenberghe, \emph{Convex optimization}, Cambridge University
  Press, Cambridge, 2004.

\bibitem[CP01]{Cauwenberghs_01_1}
G.~Cauwenberghs and T.~Poggio, \emph{Incremental and decremental support vector
  machine learning}, Advances in Neural Information Processing Systems,
  vol.~13, MIT Press, 2001, pp.~409--415.

\bibitem[CK00]{Colonius_00_1}
F.~Colonius and W.~Kliemann, \emph{The Dynamics of Control}, Birkh\"{a}user, 2000.

\bibitem[CV95]{Cortes_95_1}
C.~Cortes and V.~Vapnik, \emph{Support-vector networks}, Machine Learning
  \textbf{20} (1995), no.~3, 273--297.

\bibitem[CYLW98]{Clarke_98_1}
F.H. Clarke, R.J.~Stern Y.S.~Ledyaev, and P.R. Wolenski, \emph{Nonsmooth
  analysis and control theory}, Graduate Texts in Mathematics, vol. 178,
  Springer, New York, 1998.

\bibitem[Dei92]{Deimling_92_1}
K.~Deimling, \emph{Multivalued {D}ifferential {E}quations}, de Gruyter Series
  in Nonlinear Analysis and Applications, vol.~1, Walter de Gruyter \& Co.,
  Berlin, 1992.

\bibitem[DF89]{Dontchev_89_1}
A.L. Dontchev and E.M. Farkhi, \emph{Error estimates for discretized
  differential inclusion}, Computing \textbf{41} (1989), no.~4, 349--358.

\bibitem[DHV00]{Dontchev_00_1}
A.L. Dontchev, W.W. Hager, and V.M. Veliov, \emph{Second-order {R}unge--{K}utta
  approximations in control constrained optimal control}, SIAM Journal on
  Numerical Analysis \textbf{38} (2000), no.~1, 202--226.

\bibitem[DS02]{DeCoste_02_1}
D.~DeCoste and B.~Sch\"{o}lkopf, \emph{Training invariant support vector
  machines}, Machine Learning \textbf{46} (2002), 161--190.

\bibitem[GHHL12]{Gerdts_12}
M.~Gerdts, R.~Henrion, D.~H\"{o}mberg and C.~Landry,
\emph{Path planning and collision avoidance for robots},
Numerical Algebra, Control and Optimization \textbf{2} (2012), no.~3, 437--463.

\bibitem[H\'{a}j08]{Hajek_08_1}
O.~H\'{a}jek, \emph{Control theory in the plane}, second ed., Lecture Notes in
  Control and Information Sciences, vol. 153, Springer, Berlin, 2008.

\bibitem[HY06]{Homburg_06_1}
A.J.~Homburg and T.~Young, \emph{Hard bifurcations in dynamical systems with bounded random perturbations}, Regular \& Chaotic Dynamics \textbf{11} (2006), no.~2, 247--258.

\bibitem[Joa97]{Joachims_97_1}
T.~Joachims, \emph{Text categorisation with support vector machines}, Technical
  report, LS VIII No. 23, University of Dortmund, 1997.

\bibitem[Joa98]{Joachims_98_1}
\bysame,
\emph{Making large-scale support vector machine learning practical}, in Advances in Kernel Methods - Support Vector Learning, Cambridge MA: MIT Press, Sch\"{o}lkopf, Burges and Smola, Eds., 1998, 169--184.


\bibitem[LRR]{Lamb_Unpub_1}
J.S.W.~Lamb, M.~Rasmussen and C.S.~Rodrigues, \emph{Topological bifurcations of minimal invariant sets for set-valued dynamical systems}, to appear in: Proceedings of the American Mathematical Society.

\bibitem[Mic86]{Micchelli_86_1}
C.A. Micchelli, \emph{Algebraic aspects of interpolation}, Approximation theory
  ({N}ew {O}rleans, {L}a., 1986), Proceedings of Symposia in Applied
  Mathematics, vol.~36, American Mathematical Society, Providence, RI, 1986,
  pp.~81--102.

\bibitem[Min10]{Minh_10_1}
H.Q. Minh, \emph{Some properties of {G}aussian reproducing kernel {H}ilbert
  spaces and their implications for function approximation and learning
  theory}, Constructive Approximation \textbf{32} (2010), no.~2, 307--338.

\bibitem[NW99]{Nocedal_99_1}
J.~Nocedal and S.J. Wright, \emph{Numerical optimization}, Springer Series in
  Operations Research, Springer, New York, 1999.

\bibitem[Pla98]{Platt_98_1}
J.C.~Platt, \emph{Fast training of support vector machines using sequential minimum optimization}, in Advances in Kernel Methods - Support Vector Learning, Cambridge MA: MIT Press, Sch\"{o}lkopf, Burges and Smola, Eds., 1998, 185--208.


\bibitem[Rie13]{Rieger_13_1}
J.~Rieger, \emph{Non-convex systems of sets for numerical analysis}, Computing
  \textbf{95} (2013), no.~1, suppl., S617--S638.

\bibitem[Rie14]{Rieger_14_1}
\bysame, \emph{Semi-implicit {E}uler schemes for ordinary differential inclusions},
SIAM J.~Numer.~Anal.~\textbf{52} (2014), no.~2, 895--914.

\bibitem[Rie]{Rieger_Unpub_2}
\bysame, \emph{Robust boundary tracking for reachable sets of nonlinear
  differential inclusion}, to appear in: Foundations of Computational
  Mathematics, DOI: 10.1007/s10208-014-9218-8.

\bibitem[SHS06]{Steinwart_06_1}
I.~Steinwart, D.~Hush, and C.~Scovel, \emph{An explicit description of the
  reproducing kernel {H}ilbert spaces of {G}aussian {RBF} kernels}, IEEE
  Transactions on Information Theory \textbf{52} (2006), no.~10, 4635--4643.

\bibitem[SS01]{Schoelkopf_01_1}
B.~Sch\"{o}lkopf and A.~Smola, \emph{Learning with kernels: Support vector
  machines, regularization, optimization and beyond}, MIT Press, 2001.

\bibitem[Ste01]{Steinwart_01_1}
I.~Steinwart, \emph{On the influence of the kernel on the consistency of
  support vector machines}, Journal of Machine Learning Research \textbf{2}
  (2001), no.~67--93.

\bibitem[TC01]{Tong_01_1}
S.~Tong and E.~Chang, \emph{Support vector machine active learning for image
  retrieval}, Proceedings of the Ninth ACM International Conference on
  Multimedia, ACM, 2001, pp.~107--118.

\bibitem[V95]{Vapnik_95_1}
V.~Vapnik,
\emph{The Nature of Statistical Learning Theory},
New York: Springer-Verlag, 1998.

\end{thebibliography}
\end{document}